# A New Challenging Curve Fitting Benchmark Test Set for Global Optimization

Peicong Cheng[1,2], Peicheng Cheng[3]

1. Dept. of Electrical Engineering & Computer Science, University of Missouri, USA;
2. School of Computer and Information Technology, Beijing Jiaotong University, China;
3. Organo Corporation, Tokyo, Japan

Corresponding Email Address: pcgkg@umsystem.edu

**Abstract**

Benchmark sets are extremely important for evaluating and developing global optimization algorithms and related solvers. A new test set named PCC benchmark is proposed especially for optimization problems of nonlinear curve fitting for the first time, with the aspiration of helping developers to investigate and compare the performance of different global optimization solvers, as well as more effective optimization algorithms could be developed. Compared with the well-known classical nonlinear curve fitting benchmark set given by the National Institute of Standards and Technology (NIST) of USA, the most distinguishable features of the PCC benchmark are small problem dimensions, unconstrained with free search domain and high level of difficulty for obtaining global optimization solutions, which make the PCC benchmark be not only suitable for validating the effectiveness of different global optimization algorithms, but also more ideal for verifying and comparing various related solvers. Seven of the world's leading global optimization solvers, including Baron, Antigone, Couenne, Lingo, Scip, Matlab-GA and 1stOpt, are employed to test NIST and PCC benchmark thoroughly in terms of both effectiveness and efficiency. The results showed that the NIST benchmark is relatively simple and not suitable for global optimization testing, meanwhile the PCC benchmark is a unique, challenging and effective test dataset for global optimization.

**Keywords** Benchmark · Nonlinear Curve Fitting · Global Optimization · Mathematical Solver

# 1  Introduction

Curve fitting is a very common mathematical method used in the area of numerical science computing, whether it is for engineering simulation modelling, socio-economic operation prediction, data mining or artificial intelligence, curve fitting plays a very significant role in many research fields.

Fundamentally, curve fitting is an optimization problem that can be generally classified into two categories according to the nature of the fitting function: linear and nonlinear. Both are solved based on least square method by applying relevant algorithms. For the linear fitting type, the processing is relatively straightforward, since linear fitting problem usually can be directly solved by matrix operations. However, handling nonlinear fitting types is significantly more intricate for linear fitting type. As shown in Equation (1), the objective function of the least squares model is to minimize the sum of squared errors (SSE) between actual and calculated dependent variable. Generally, it is necessary to employ an optimization algorithm for iterative computation procedure to obtain reasonable results, therefore, the optimization algorithm becomes the key for solving nonlinear fitting problems.

$$\text{Min. } SSE = \sum_{i=1}^{n}(y_i - \dot{y}_i)^2 \qquad (1)$$

Where $y$ and $\dot{y}$ are the actual and calculated dependent variable values, respectively, and SSE is the sum of squared error, its theoretical lower bound value or the best ideal result is zero.

The nonlinear curve fitting problem is actually treated as an optimization problem, and thus the selection and application of the optimization algorithm have a huge impact on the quality of the fitting results, while on the other hand, the curve fitting problem can also be used to examine and assess the goodness of the optimization algorithms or relevant solvers.

Nonlinear optimization algorithms usually can be divided into two types: local optimization algorithm (LOA) and global optimization algorithm (GOA). LOA, as the name implied, can get the local solution usually. The typical LOAs including Newton algorithm, Quasi-Newtonian algorithm, conjugate gradient algorithm, steepest descent method, etc. The LOAs are characterized by high efficiency but rely heavily on the guessing of the initial start-values of the parameters, the quality of initial start-values guessed has a decisive effect on final convergence of the iterative computation, this could even be thought as the biggest shortcoming of LOA. While the GOA, as opposed to the LOA, theoretically has the most important feature of no longer relying on guessing of initial start-values but achieving the global optimization solution. Generally, there are two types of the GOA: deterministic and stochastic algorithms, the former including such Branch-and-bound (BB) algorithms [6, 23], DIRECT-type algorithms [11] and multi-start (MS) search algorithms [8]; while the later including genetic algorithms [10], Tabu Search [5] and simulated annealing [13], as well as a large number of heuristic optimization algorithms such as particle swarm algorithms [9, 22] and ant colony algorithms [3, 16]. Theoretically, for nonlinear curve fitting problems, except the drawback of relative low efficiency, GOA indeed have huge advantages over LOA in terms of ease of use and computation effect. However, the truth is, rather than GOA, the traditional LOA, for example, the Levenberg-Marquardt algorithm, are still the major default algorithms in the fields of nonlinear curve fitting. This is mainly because GOA does not guarantee to obtain the global optimal solution for all problems in practice, although it is theoretically possible.

Although there are already some excellent global optimization solvers available on the market today, there are rarely any comprehensive and universal tests and comparison studies on the effectiveness and efficiency of these professional global optimization solvers in the field of nonlinear curve fitting. One of the main reasons is the lack of a high standard benchmark test set for curve fitting problems. Though there are some benchmarks for nonlinear fitting problems, such as the classical curve fitting problems from National Institute of Standards and Technology (NIST) [21], they were proposed one or several decades ago and mainly used to test LOA, and are relatively simple and have no practical significance for GOA. Therefore, a new up-to-date nonlinear curve fitting benchmark for testing and evaluating global optimization solvers is urgently needed today, which is the primary reason that PCC benchmark is proposed this time. The PCC benchmark will be not only useful for testing and evaluating the performance of existing global solvers, but also helpful for researchers to develop more efficient global optimization algorithms.

The PCC benchmark, summarized and proposed by the authors in this paper, contains 38 fitting problems with a high degree of difficulty. In addition to NIST benchmark, seven well-known optimization solvers, Baron, Antigone, Couenne, Lingo, Scip, Matlab-GA and 1stOpt, have been employed to examine these two benchmarks. Two objectives are hoped to find out by this study, one is whether the proposed benchmark is challenging for current leading global optimization solvers, another is how the performance of those solvers will be on these two fitting problems benchmark.

## 2   Benchmark of Nonlinear Curve Fitting

## 2.1 NIST Benchmark

NIST nonlinear curve fitting benchmark is one of the most well-known benchmarks with a set of 27 nonlinear curve fitting problems as shown in Table 1. The benchmark originally targets on evaluating all sorts of LOAs and provides two sets of initial start-values for each problem as well as the certified correct results, aiming to check whether the given LOAs can finally obtain the correct results by using either of the two initial start-values provided and assess the performance of the algorithms based on the results. The benchmark problems are divided into three groups according to their difficulty: high, average and lower. In the field of curve fitting, almost all famous data analysis packages in the world take this benchmark as the standard for verifying and judging the quality of their optimization algorithms. Although the original intention of this benchmark is to test LOA, GOA can also apply on the NIST benchmark and evaluate the performance, since GOA theoretically has the characteristics and advantages of global convergence without depending on the initial start-values. Therefore, it is an interesting thing to discover the performance of GOA from the selected seven solvers on this benchmark. Without using the initial start-values provided by NIST, can these solvers obtain the certified solutions without too much effort? Is it possible to employ GOA as a superior substitute of the LOA in real practice? Such studies rarely have been done before, so it is very practical and necessary to conduct a comparative study in this issue. The NIST benchmark can be referred to [21] for details.

**Table 1** NIST Benchmark

| Dataset Name | Level of Difficulty | Model Classification | No. of Parameters | No. of Observations | Source |
|---|---|---|---|---|---|
| Misra1a | Lower | Exponential | 2 | 14 | Observed |
| Chwirut2 | Lower | Exponential | 3 | 54 | Observed |
| Chwirut1 | Lower | Exponential | 3 | 214 | Observed |
| Lanczos3 | Lower | Exponential | 6 | 24 | Generated |
| Gauss1 | Lower | Exponential | 8 | 250 | Generated |
| Gauss2 | Lower | Exponential | 8 | 250 | Generated |
| DanWood | Lower | Miscellaneous | 2 | 6 | Observed |
| Misra1b | Lower | Miscellaneous | 2 | 14 | Observed |
| Kirby2 | Average | Rational | 5 | 151 | Observed |
| Hahn1 | Average | Rational | 7 | 236 | Observed |
| Nelson | Average | Exponential | 3 | 128 | Observed |
| MGH17 | Average | Exponential | 5 | 33 | Generated |
| Lanczos1 | Average | Exponential | 6 | 24 | Generated |
| Lanczos2 | Average | Exponential | 6 | 24 | Generated |
| Gauss3 | Average | Exponential | 8 | 250 | Generated |
| Misra1c | Average | Miscellaneous | 2 | 14 | Observed |
| Misra1d | Average | Miscellaneous | 2 | 14 | Observed |
| Roszman1 | Average | Miscellaneous | 4 | 25 | Observed |
| ENSO | Average | Miscellaneous | 9 | 168 | Observed |
| MGH09 | Higher | Rational | 4 | 11 | Generated |
| Thurber | Higher | Rational | 7 | 37 | Observed |
| BoxBod | Higher | Exponential | 2 | 6 | Observed |
| Rat42 | Higher | Exponential | 3 | 9 | Observed |
| MGH10 | Higher | Exponential | 3 | 16 | Generated |
| Eckerle4 | Higher | Exponential | 3 | 35 | Observed |
| Rat43 | Higher | Exponential | 4 | 15 | Observed |
| Bennett5 | Higher | Miscellaneous | 3 | 154 | Observed |

## 2.2 PCC Benchmark

The NIST benchmark test set enjoys high reputation in both popularity and recognition, but clearly suffers from the following two problems: one is that the original intention of this benchmark is to be proposed specifically for testing the LOA, which implies it may not be appropriate enough for testing GOA; another is that this benchmark had been proposed more than thirty years and was indeed difficult and representative test set for developing and evaluating the curve fitting algorithms at that time. However, with the rapid development and progress of modern science and technology, is this benchmark still challenging? Is it relatively simple for today's GOA solvers? If the certified solutions can be obtained easily by GOA solvers, such benchmark will then be unconvincing and meaningless for test and verify GOA solvers nowadays.

Currently, there are lots of benchmark test datasets for global optimization such as the well-known GLOBALLib [18], MINLPLib [20] and DIRECTGOLib [25]. Unfortunately, none of them are specifically for checking and testing the performance of GOA on nonlinear curve fitting. Moreover, most of those problems typically exhibit two characteristics: the first is that the search domains are provided, which result in the reduction of difficulties of those problems are reduced greatly; The second is that the function types of part of those problems are most-like man-made intentionally, i.e., they are not derived from reality, but rather created for the sake of the problem. For some problems, the optimal solution can be directly solved just by simple observation, but turn out to be tough for GOA solvers. In view of this, the authors have collected, summarized and organized 38 nonlinear curve fitting problems from practical experience, named PCC Benchmark test dataset, particularly for nonlinear curve fitting. All problems contain only one independent and one dependent variable, the fitting functions are not complicated, involving only primary math functions such as exponential, logarithmic and power functions, the problem dimensions or unknown parameter numbers varies from 2 to 8, all problems are low-dimensional and unconstrained with free search domain, the global optimization solution for each problem is hard to be identified. Moreover, unlike traditional optimization problems, such as the well-known Schwefel and Shubert functions represented in Eq. (2) and Eq. (3), which usually have impressed and straightforward 3D charts, as shown in Fig. 1 and Fig. 2, that clearly showing multiple local optimal solution, PCC benchmark problems are difficult to visually represent the global and local optimal solutions of each problem through the same chart because of the free search domain and the properties of fitting functions.

- Schwefel Function, the search domain is [-500, 500].

$$Min. \sum_{i=1}^{n}\left(-x_i \sin\left(\sqrt{|x_i|}\right)\right) \quad (2)$$

- Shubert Function, the search domain is [-10, 10].

$$Min. \left(\sum_{i=1}^{2}\left(i\cos((i+1)x_1+i)\right)\right)\left(\sum_{i=1}^{2}\left(i\cos((i+1)x_2+i)\right)\right) \quad (3)$$

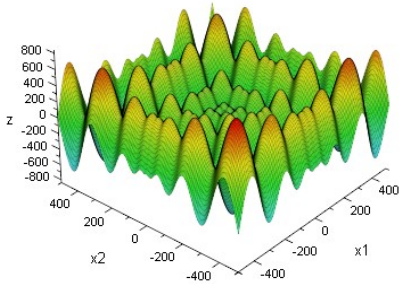 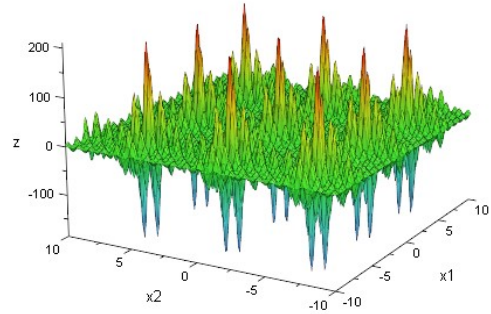

Fig.1 Surface Chart of Schwefel Function    Fig.2 Surface Chart of Shubert Function

The PCC benchmark is published for the first time, there's no certified global solution for each problem, the best optimal solution obtained so far for each problem will be considered as the global solution.

The fitting function and corresponded data for all 38 test problems are shown in Table 2, where x and y represent the independent and dependent variables, respectively, and b is the parameter vector.

**Table 2** PCC Benchmark

| No. | Fitting function and data | No. of Parameters |
|---|---|---|
| 1 | $y = b_1 x^{b_2}$<br>x=[20,25,27,30,35,38.5]; y=[4.62e-8,1.02e-7,1.27e-7,2.14e-7,3.56e-7,4.69e-7]; | 2 |
| 2 | $y = b_1 - b_2 \exp\left(b_3 x + \dfrac{1}{b_4 x^{b_5} + b_6 + 1}\right)$<br>x=[0.9,1.5,13.8,19.8,24.1,28.2,35.2,60.3,74.6,81.3];<br>y=[455.2,428.6,124.1,67.3,43.2,28.1,13.1,-0.4,-1.3,-1.5]; | 6 |
| 3 | $$y = \frac{b_1}{1 + b_2\exp(-b_3 x - b_4 x^2)} + \frac{2b_4}{(1 + b_5\exp(-b_6 x))^{b_4}}$$ x = [-6.7,-6.3,-5.9,-5.5,-5.1,-4.7,-4.3,-3.7,-3.3,-2.3,-1.9,-1.5,-1.2,-0.7,-0.2,0.3,0.7]; y = [0.06,0.1,0.14,0.22,0.34,0.5,0.62,0.72,0.75,0.87,0.95,1.05,1.15,1.35,1.45,1.5,1.54]; | 6 |
| 4 | $$y = p_1 p_2^x (x + p_3)^{p_4} + \frac{p_5 x}{1 + p_6 x}$$ x = [0,2,4,6,8,10]; y = [0.0242,0.2792,0.1386,0.0238,0.0034,0.002]; | 6 |
| 5 | $$y = \frac{p_1}{1 + p_2\exp(-p_3 x) + p_4 x^{p_5}} + p_6 x$$ x = [1,2,3,4,5,6,7,8,9,10]; y = [161.68,170.99,179.42,185.64,189.29,192.58,196.32,205.40,240.02,293.98]; | 6 |
| 6 | $$y = b_1\exp\left(\frac{b_2}{x} - b_3 x + b_4 x^{b_5}\right) + b_6 x$$ x = [0.17,0.33,0.83,1.67,2.5,3.33,5]; y = [6.62,5.58,4.63,5.08,5.42,5.82,6.78]; | 6 |
| 7 | $$y = p_1 + p_2\exp(p_3\exp(p_4 x) + p_5 x) + p_6 x$$ x = [20,21,22,23,24,25,26,27,28,29,30]; y = [56.7,59.8,64.8,66.7,66.2,65.2,64.2,63.5,63.0,62.7,62.5] | 6 |
| 8 | $$y = b_1\exp(b_2 x) + b_3\exp\left(\frac{b_4}{x}\right) + b_5$$ x = [250,275,300,325,350]; y = [2,6,15,19,37]; | 5 |
| 9 | $$y = \frac{b_3}{1 + \exp(b_1 x + b_5)} + \frac{b_3}{1 + b_4\exp(b_2 x)} + b_6 x + b_5$$ x=[1.75,2.25,2.493,2.5,2.563,2.686,2.885]; y=[260,320,381,440,460,500,570]; | 6 |
| 10 | $$y = b_1 + b_2\exp(b_3(x - b_4)^2) + b_5\exp(b_6(x - b_7)^2)$$ x=[ 1,2,4,6,8,10,15,20,30,40,60,90,120,150,180]; y=[ 0,0,0,0,0.20640,0.06304,0.08540,0.56381,0.94651,1.46233,1.73185,1.70072,1.82921,2.09236,2.26539]; | 7 |
| 11 | $$y = \frac{b_1 + b_2 x + b_3\exp(b_4 x)}{1 + b_5 x + b_6\exp(b_7 x)}$$ x=[0,10,20,30,40,50,60,70,80,90]; y=[0,4.808,18.367,37.920,58.726,74.756,81.454,79.025,71.244,61.629]; | 7 |
| 12 | $$y = b_1 + b_2 x + b_3\exp\left(-(b_4(x - b_5))^2\right) + b_6\exp\left(-(b_7(x - b_8))^2\right)$$ x=[0.975,0.985,0.995,1.015,1.025,1.03,1.035,1.04,1.045,1.05,1.055,1.06,1.07,1.09,1.1,1.13]; y=[306,323,438,828,1068,1195,1267,1257,1243,1246,1235,1126,889,474,327,237]; | 8 |
| 13 | $$y = b_1 + b_2\left(\frac{b_3}{1 + \exp(b_5(x - b_4))} + \frac{1 - b_3}{1 + \exp(b_4(x - b_5))}\right)$$ x=[-14,-12.45,-9.36,-6.26,-3.17,0.07,3.03,6.12,9.22,12.31,14]; y=[166,123,47.5,-8.85,-30.6,0.842,30.6,10.6,-44.5,-116,-165]; | 5 |
| 14 | $$y = b_1 + b_2 x^{b_3} + b_4 x^{b_5}$$ | 5 |
| | x=[5,10,15,20,25,30,35,40,45,50,55,60,65,70,75,80,85,90,95]; <br> y=[3.3,6.5,9.2,11.9,14.5,17.0,19.4,21.7,23.9,25.9,27.6,29.2,30.0,30.3,30.0,29.2,26.0,21.0,12.0]; | |
| 15 | $$y = \frac{b_1}{1 + b_2 \exp(b_3 x)} + \frac{b_4}{1 + b_5 \exp(b_6 x)}$$ <br> x=[ 0,10,20,30,40,50,60,70,80,90,100]; y=[ 92,93,94.5,96,98,100,101,100,93,83,74]; | 6 |
| 16 | $$y = b_1 \exp\left(\frac{b_2}{x} - b_3 x + b_4 x^{b_5}\right) + b_6 x$$ <br> x=[ 0.17,0.33,0.83,1.67,2.5,3.33,5]; y=[6.62,5.58,4.63,5.08,5.42,5.82,6.78]; | 6 |
| 17 | $$y = b_1 \left(\left(\frac{x - b_2}{(1 - b_2)x + b_3}\right)^2 + b_4 x + 2.1\right)^{b_5}$$ <br> x=[1,0.963,0.940,0.810,0.639,0.564,0.495,0.456,0.374,0.324,0.289,0.248]; <br> y=[0,7.642,84.369,730.794,943.654,1203.544,1489.304,1814.872,2541.681,3399.595,4385.219,5469.341]; | 5 |
| 18 | $$y = b_1 \exp(b_2 x) + b_3 \exp\left(b_4 x^{0.5} + \frac{b_5}{x}\right) + b_6$$ <br> x=[0.392,0.793,1.189,2.39,4.787,7.191,9.385,11.784]; <br> y=[1.63,6.79,6.79,3.88,1.83,1.25,0.86,0.624]; | 6 |
| 19 | $$y = \frac{b_1}{1 + b_2 x + b_3 x^2} + b_4 \exp(b_5 x)$$ <br> x=[1.5,2.5,3.5,4.5,5.5,6.5]; y=[1.69,1.636,1.387,1.254,1.444,0.1631]; | 5 |
| 20 | $$y = b_1 \exp(-b_2 \exp(b_3 x)) + b_4$$ <br> x=[11.95,12.22,12.88,13.50,13.96,14.25,14.32,14.52,14.78,15.04,15.30]; <br> y=[7.66,7.66,7.52,7.10,6.50,6.01,5.90,5.61,5.37,5.28,5.25]; | 4 |
| 21 | $$y = b_1 + b_2 \exp(-\exp(0.03(x - b_3)^{b_4}) - b_5(x - b_3) + 1) + b_6 x$$ <br> x=[0,0.017,0.034,0.068,0.102,0.136,0.203,0.237,0.305,0.373,0.441,0.492]; <br> y=[3.33,5,7.4,8.61,7.72,5.98,3.16,2.4,1.15,0.56,0.33,0.26]; | 6 |
| 22 | $$y = \frac{b_1 + b_2 x + b_3 \exp(b_4 x)}{1 + b_5 \exp(b_6 x)}$$ <br> x=[22,52,80,109,125,142,172,201,230,261,295]; <br> y=[220.7,210.1,187,200,310,318,335,325,329,328,325]; | 6 |
| 23 | $$y = \frac{b_1}{1 + b_2 \exp(-b_3 x) + b_4 x} + b_5 x^{b_6}$$ <br> x=[2,3,4,5,6,7,8,9,10,11,12,13]; y=[0.3,2.1,7.7,13.8,12.9,12.5,9.8,9,9.1,9.7,9.2,10.8]; | 6 |
| 24 | $$y = b_1 + b_2 x + b_3 x^{b_4} + \frac{b_5}{b_2 + x}$$ <br> x=[140,150,180,220,260,280,290,298.15,360,600,1000]; <br> y=[118.73,125.55,144.14,164.31,181.26,188.78,192.55,195.18,213.97,252.4,288.67]; | 5 |
| 25 | $$y = b_1 + \frac{b_2}{1 + b_4 \exp(b_3 x) + b_5 \exp(\exp(b_6 x))}$$ | 6 |
| | x=[-7,-6.5,-6,-5.5,-5,-4.5,-4,-3.5,-3.0]; y=[1,5,14,19,23,26,26.1,27,27.1]; | |
| 26 | $$y = b_1 \ln(x + b_2) + b_3 \exp(-b_4 x) + b_5 \exp(-b_6 x)$$ x=[ 300,600,900,1200,1800,3600,5400,7200,27000]; y=[ 39.059,23.585,18.418,15.825,13.208,10.363,8.812,7.285,2.007]; | 6 |
| 27 | $$y = b_1 + b_2 \exp\left(-0.5\big(b_3(x - b_4)\big)^2\right) - b_2 \exp\left(-0.5\big(2b_3(x - b_5)\big)^2\right)$$ x=[0,54,174,250,282,330,355]; y=[98,17.95,46.51,48.37,11.16,74.4,100]; | 5 |
| 28 | $$y = b_1 \exp(b_2(x - b_3)^2) + b_4 \exp(b_5(x - b_6)^2)$$ x=[2015,2016,2017,2018,2019,2021,2023,2025,2027,2028,2029,2030,2031,2033,2035,2037,2039,2040]; y=[23.268,24.46,24.873,24.959,24.925,24.775,24.621,24.484,24.374,24.335,24.311,24.304,24.318,24.417,24.621,24.931,25.338,25.571]; | 6 |
| 29 | $$y = b_1 \exp\big(-\exp(-b_2(x - b_3)) - b_1\big) + b_4 \exp(b_5 x - b_3)$$ x=[12.64,13.344,14.048,14.752,15.104,15.2,15.232,15.456,15.808,16.864,17.76]; y=[-0.040,-0.033,-0.035,-0.031,0.076,0.406,0.624,0.965,1.001,0.998,0.997]; | 5 |
| 30 | $$y = \frac{b_1 + b_2 \exp(b_3 x)}{1 + b_4 x + b_5 \exp(b_6 x)}$$ x=[23.5,64.8,145.3,185.2,267.1,389.7,513.7,638.2,680.3,721.4,763.3,887.5,1013,1138.6,1263.8,1305.8]; y=[55,122,187,190,173,138,111,92,86,82,78,68,60,53,48,47]; | 6 |
| 31 | $$y = \frac{1 + b_1 x^{b_2}}{b_3 + \exp(b_4 + b_5 x + b_6 x^{b_7})}$$ x=[0.1001,0.1602,0.2288,0.3079,0.4003,0.5094,0.6403,0.8002]; y=[87.7,89.8,100.2,121.2,131.3,138.9,135.4,121.4]; | 7 |
| 32 | $$y = b_1 + \exp(-b_1 - b_2 b_3 x) + b_3 x^{b_4 + b_5 x}$$ x=[0.25,0.5,0.75,1,1.5,2.5,3.5,4.5,5,6,8,10,12,14,16]; y=[15,34,37.5,41,41,34,29,25,20.5,19,14,9,6,3.5,2]; | 5 |
| 33 | $$y = \left(b_1 b_2 - b_3 \exp(b_1 - b_4 x) + \frac{b_2 b_3}{x}\right)^{b_5}$$ x=[1,2,4,6,8]; y=[100,140,160,170,175]; | 5 |
| 34 | $$y = \frac{b_1 b_2 + b_3 \exp(b_4 x)}{b_2 + \exp(b_4 x) + \frac{b_5}{x + b_1}}$$ x=[0,10,20,30,40,50,60,70,80,90]; y=[22.260,21.526,20.589,19.092,17.312,15.425,13.854,12.541,11.697,11.294]; | 5 |
| 35 | $$y = b_1 + b_2 \exp(b_3 x + b_4 x^{b_5})$$ x=[204.11,219.64,243.34,281.01,312.05,331.71,341.49,349.16,379.08,420.91,430.99]; y=[64.33,69.17,71.18,78.55,82.57,81.23,77.37,73.13,60.59,23.04,2.98]; | 5 |
| 36 | $$y = b_1 \exp\left(\frac{b_2}{(b_3 + x)^2} - b_4 x\right) + b_5$$ x=[298.15,303.15,308.15,313.15,315.5,318.15,323.15,328.15,333.15,338.15]; y=[114.6,89.68,68.43,50.21,3.95,40.09,33.42,27.01,21.59,17.53]; | 5 |

| 37 | $y = b_1 + b_2\exp(b_3 x) + b_4\exp\left(\dfrac{b_5}{x}\right)$ | 5 |
|---|---|---|
| | x=[73.1,94,100,130,150,184,192,200]; y=[0.25,0.6,0.65,0.9,0.91,0.94,0.95,1]; | |
| 38 | $y = b_1 + \dfrac{b_2}{1 + b_3 x + b_4\exp(b_5 x)}$ | 5 |
| | x=[-24536,-26603,-27424,-28573,-29204,-29428,-29822,-30065]; y=[0,107,216,414,643,785,1044,1285]; | |

# 3 Global Optimization Solvers

The seven global optimization solvers used in this study are briefly described below.

The Baron solver was originally developed by DR. Sahinidis, an academician of the American Academy of Engineering and current professor at the Georgia Institute of Technology Carnegie Mellon University, in 1996 when he worked at the University of Illinois at Urbana-Champaign, and after more than 20 years of continuous improvement and upgrading, it has been recognized as the most advanced and powerful global optimization solver today [7, 14, 19, 12].

The Antigone solver is an optimization solver jointly developed by Prof. Floudas from Princeton University of U.S.A and Prof. Misener from Imperial College of UK, and first published in 2013, it focuses on global solutions for continuous and integer programming, and enjoys a high reputation in the market [17].

Couenne solver is a nonlinear global optimization solver developed by COIN-OR in the United States, originated from a collaborative project between IBM and Carnegie Mellon University in 2006, the solver is open source with rich API for easy to use [1].

As an optimization software package developed by Lindo Company in 1981, Lingo software has its own model language, it is easy to get start with and widely used, and is chosen by more than half of the world's top 500 companies. Lingo is also the high-recommended software tool for various Mathematical modeling competition [4].

SCIP originated from ZIB Research Center in Berlin, Germany, and has evolved from being mainly used for solving mixed integer programming problems to solving all kinds of optimization problems, including nonlinear global optimization problems, and is a highly evaluated open-source software with rich interfaces and support for C/C++, Java, Python, MATLAB and other programming environments [2, 24].

Matlab is a commercial mathematical software produced by MathWorks Corporation in the United States. Together with Mathematica and Maple, it is known as the three major mathematical software today. Matlab is also recognized as the best solver in the areas of numerical computation, it involves and integrates a number of toolboxes for various purposes, and the genetic algorithm toolbox (Matlab-GA) is one of its built-in toolboxes [15].

1stOpt is an optimization solver developed by 7D-Soft High Technology Inc. from China, the most impressive feature of 1stOpt is its unique Universal Global Optimization (UGO) algorithm, which greatly improves the ability of global optimization; moreover, it's simple user interface and easy-to-use model language attracted many users in a short period of time, and there are about 10,000 scientific research papers have used 1stOpt as their data processing and analysis tool [26].

GAMS is a general-purpose optimization platform that integrates many well-known optimization solvers and provides a unified user interface and model language for all different solvers. One advantage of GAMS is that users simply need to master GAMS model language to solve problem in different supportive solvers without having to understand specific language of each solver. In another word, users are capable to use same problem code in different solvers, with only one line of code need to be modified to select different solvers, which reduces the application cost and greatly improves user's efficiency.

The first five solvers mentioned above are already included in the GAMS platform, so their computation will be implemented on GAMS platform in this study. 1stOpt and Matlab will use their own platform. The trial version of GAMS 43.4.1 is employed this time, which can be downloaded freely from the GAMS official website. It should be noted that this version of GAMS does not include the Couenne solver by default, it is essential to separately download the relevant package and configure the settings before utilizing the Couenne solver seamlessly on GAMS platform. The versions of Matlab and 1stOpt are 2020B and 10.0, respectively.

## 4 Test Comparison

Based on the above NIST and PCC benchmark test datasets, each solver was conducted for computation and comparison in the same hardware environment. Because the first five solvers in GAMS use deterministic algorithms and the computational outputs are stable and unchanging, only 1stOpt and Matlab-GA were run independently 10 times for each test problem, and their respective best outputs were taken as their final results. Furthermore, since each solver comes with numerous configuration options and different options can have an impact on the final results, considering the fact that default configurations are generally suitable for most cases and have the highest expectation of getting the best result, for the convenience in this study, the default configuration options are adopted for each solver. For Matlab-GA, the population size is taken as 2000, and the maximum number of iterations is set to 3000.

### 4.1 Test Results of NIST Benchmark

Certificated correct results have been provided for each problem in the NIST benchmark, which making it convenient to verify the effect of each solver. The initial parameter values provided by NIST are completely disregarded. Instead, the initial parameter values are automatically determined by each solver based on their algorithms, without manual settings. The calculation results are shown in Tables 3, 4, and 5, respectively. "✓" indicates the calculation result matches the certificated one, while "✘" means opposite. There's one more thing that should be noted in advance, which is for the Baron and Antigone, the reason they failed in some cases results from their lack of supporting math functions such as trigonometric functions of "sin" and "cos". The problems of "Roszman1" and "ENSO" in the NIST benchmark are in such situation.

For the lower difficulty dataset, 1stOpt performed the best with 100% correct solutions, while Baron, Lingo, Antigone, Couenne and Scip could achieve with the average success rate of 87.5%, 62.5%, 50%, 50% and 37.5% respectively. Matlab-GA exhibited less favorable performance compared to other solvers, achieving a success rate of 12.5%, with only one problem solved approximately. For the average difficulty dataset, 1stOpt was in 100% success rate, followed by Baron with 72.72%, 63.64% for Couenne, 54.54% for both Lingo and Scip, 45.45% for Antigone and finally 9.09% for Matlab-GA. For the higher difficulty dataset, both Baron and 1stOpt performed the best with a 100% success, Antigone was followed with 87.5%, Lingo, Couenne and Scip were 75%, 62.5 and 50%, respectively, and Matlab-GA was the last, with 25% only. The overall success rate of all these seven solvers on NIST benchmark is 61.9%.

In terms of overall effect, see Fig. 3. 1stOpt performs the best among the seven solvers and is the only one that can find all correct answers with random initial start-values, followed by Baron, Lingo, Antigone, Couenne, Scip and Matlab-GA, corresponding to an average success rate of 86.74%, 64.01%, 60.98%, 58.71%, and 47.34% and 15.53%, respectively. For the overall efficiency, the average time to successfully solve the test problem for each solver was taken for evaluation, see Fig. 4. The shorter the time, the more efficient the solver, and the order of efficiency from highest to lowest is Couenne, Baron, 1stOpt, Antigone, Scip, Matlab-GA, and Lingo.

In addition, through the original benchmark dataset has been divided into three groups based on their difficulty, from the actual performance results, it seems that the difficulty index is only suitable for testing the LOA, but not for GOA solvers. To be specific, in the higher difficulty group, its overall average success rate is 71.43%, which is higher than the 57.14% of the lower difficulty group. This is mainly because the initial values given for LOAs in each problem have also been taken into account when the developer of NIST benchmark labeled the difficulty index for each problem. For GOAs, without using the given initial value, problems' difficulty may change dramatically.

**Table 3** Test results of lower difficulty problems

| Problem | Solver | Baron | Antigone | Couenne | Lingo | Scip | Matlab-GA | 1stOpt |
|---|---|---|---|---|---|---|---|---|
| Misra1a | Success | ✓ | ✓ | ✗ | ✓ | ✓ | ✗ | ✓ |
|  | Time (s) | 0.20 | 0.5 | - | 120.0 | 15.0 | - | 0.5 |
| Chwirut2 | Success | ✓ | ✓ | ✓ | ✓ | ✓ | ✗ | ✓ |
|  | Time (s) | 0.5 | 5.0 | 0.5 | 0.5 | 0.5 | - | 0.6 |
| Chwirut1 | Success | ✓ | ✓ | ✓ | ✓ | ✗ | ✗ | ✓ |
|  | Time (s) | 0.5 | 5.0 | 0.5 | 0.5 | - | - | 1.5 |
| Lanczos3 | Success | ✓ | ✗ | ✓ | ✗ | ✗ | ✗ | ✓ |
|  | Time (s) | 0.5 | - | 0.5 | - | - | - | 5.0 |
| Gauss1 | Success | ✗ | ✗ | ✗ | ✗ | ✗ | ✗ | ✓ |
|  | Time (s) | - | - | - | - | - | - | 19 |
| Gauss2 | Success | ✓ | ✗ | ✗ | ✗ | ✗ | ✗ | ✓ |
|  | Time (s) | 2.66 | - | - | - | - | - | 50.0 |
| DanWood | Success | ✓ | ✓ | ✗ | ✓ | ✓ | ✓ | ✓ |
|  | Time (s) | 0.1 | 0.1 | - | 0.1 | 0.1 | 2.3 | 0.1 |
| Misra1b | Success | ✓ | ✗ | ✓ | ✓ | ✗ | ✗ | ✓ |
|  | Time (s) | 0.17 | - | 16.5 | 120.0 | - | - | 0.3 |
| Success rate (%) |  | 87.5 | 50.0 | 50.0 | 62.5 | 37.5 | 12.5 | 100.0 |
| Average time used (s) |  | 0.66 | 3.37 | 4.5 | 48.22 | 5.20 | 2.30 | 9.63 |

**Table 4** Test results of average difficulty problems

| Problem | Solver | Baron | Antigone | Couenne | Lingo | Scip | Matlab-GA | 1stOpt |
|---|---|---|---|---|---|---|---|---|
| Kirby2 | Success | ✓ | ✓ | ✓ | ✓ | ✓ | ✗ | ✓ |
|  | Time (s) | 2.0 | 8.0 | 1.5 | 1.0 | 1.5 | - | 3.5 |
| Hahn1 | Success | ✓ | ✓ | ✓ | ✓ | ✓ | ✗ | ✓ |
|  | Time (s) | 0.3 | 0.3 | 0.5 | 0.5 | 0.5 | - | 3.0 |
| Nelson | Success | ✓ | ✓ | ✓ | ✓ | ✓ | ✗ | ✓ |
|  | Time (s) | 0.2 | 50.0 | 0.4 | 60.0 | 0.3 | - | 2.0 |
| MGH17 | Success | ✓ | ✓ | ✓ | ✗ | ✓ | ✗ | ✓ |
|  | Time (s) | 1.2 | 0.3 | 0.5 | - | 0.5 | - | 2.5 |

| Problem | | Baron | Antigone | Couenne | Lingo | Scip | Matlab-GA | 1stOpt |
|---|---|---|---|---|---|---|---|---|
| Lanczos1 | Success | ✓ | ✗ | ✓ | ✗ | ✗ | ✗ | ✓ |
| | Time (s) | 0.5 | - | 0.5 | - | - | - | 1.2 |
| Lanczos2 | Success | ✓ | ✗ | ✓ | ✗ | ✗ | ✗ | ✓ |
| | Time (s) | 0.5 | - | 0.5 | - | - | - | 1.2 |
| Gauss3 | Success | ✗ | ✗ | ✗ | ✗ | ✗ | ✗ | ✓ |
| | Time (s) | - | - | - | - | - | - | 16.0 |
| Misra1c | Success | ✓ | ✓ | ✓ | ✓ | ✓ | ✗ | ✓ |
| | Time (s) | 0.5 | 2.5 | 1.5 | 1.5 | 1.5 | - | 0.5 |
| Misra1d | Success | ✓ | ✗ | ✗ | ✓ | ✓ | ✗ | ✓ |
| | Time (s) | 0.5 | - | - | 0.5 | 16.0 | - | 0.5 |
| Roszman1 | Success | ✗ | ✗ | ✗ | ✓ | ✗ | ✗ | ✓ |
| | Time (s) | - | - | - | 0.5 | - | - | 1.0 |
| ENSO | Success | ✗ | ✗ | ✗ | ✗ | ✗ | ✓ | ✓ |
| | Time (s) | - | - | - | - | - | 17.0 | 18.0 |
| Success rate (%) | | 72.72 | 45.45 | 63.64 | 54.54 | 54.54 | 9.09 | 100 |
| Average time used (s) | | 0.71 | 12.22 | 0.77 | 10.67 | 3.38 | 17.00 | 4.49 |

**Table 5** Test results of higher difficulty problems

| Problem | | Baron | Antigone | Couenne | Lingo | Scip | Matlab-GA | 1stOpt |
|---|---|---|---|---|---|---|---|---|
| MGH09 | Success | ✓ | ✓ | ✗ | ✓ | ✓ | ✗ | ✓ |
| | Time (s) | 0.5 | 0.5 | - | 0.5 | 0.5 | - | 0.5 |
| Thurber | Success | ✓ | ✓ | ✗ | ✓ | ✓ | ✗ | ✓ |
| | Time (s) | 0.5 | 0.5 | - | 69.0 | 0.5 | - | 1.5 |
| BoxBod | Success | ✓ | ✓ | ✗ | ✓ | ✓ | ✓ | ✓ |
| | Time (s) | 0.5 | 0.5 | - | 0.2 | 0.5 | 1.5 | 1.0 |
| Rat42 | Success | ✓ | ✓ | ✓ | ✓ | ✓ | ✗ | ✓ |
| | Time (s) | 0.5 | 0.5 | 0.5 | 0.5 | 0.5 | - | 0.5 |
| MGH10 | Success | ✓ | ✓ | ✓ | ✓ | ✗ | ✗ | ✓ |
| | Time (s) | 0.25 | 0.5 | 0.5 | 0.5 | - | - | 1.0 |
| Eckerle4 | Success | ✓ | ✓ | ✓ | ✗ | ✗ | ✓ | ✓ |
| | Time (s) | 2.0 | 1.0 | 1.0 | - | - | 2.8 | 1.0 |
| Rat43 | Success | ✓ | ✗ | ✓ | ✗ | ✗ | ✗ | ✓ |
| | Time (s) | 0.5 | - | 1.0 | - | - | - | 1.0 |
| Bennett5 | Success | ✓ | ✓ | ✓ | ✓ | ✗ | ✗ | ✓ |
| | Time (s) | 1.0 | 1.0 | 1.0 | 1.0 | - | - | 10.0 |
| Success rate (%) | | 100 | 87.5 | 62.5 | 75 | 50 | 25 | 100 |
| Average time used (s) | | 0.72 | 0.64 | 0.80 | 11.95 | 0.50 | 2.15 | 2.06 |

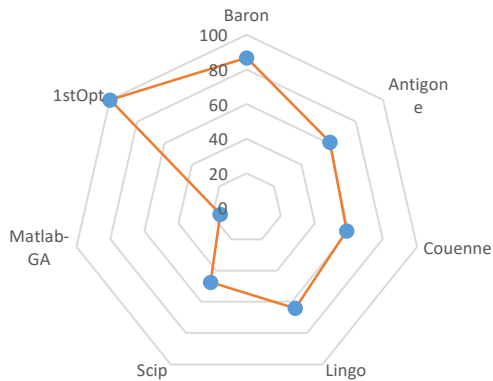
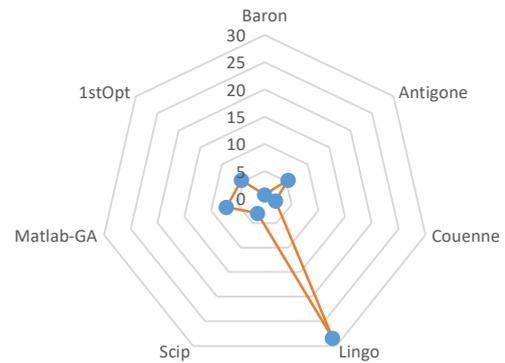

Fig. 3 Effect: Success Rate (%) of Each Solver   Fig. 4 Efficiency: Ave. Time (s) of Each Solver

## 4.2 Test Results of PCC Benchmark

In contrast to the NIST benchmark, the 38-problems PCC benchmark dataset does not know the theoretical correct answers or certified global optimal solutions, the best result obtained for each solver so far will be considered as the global optimal solution. Due to the high difficulty of global optimization problem, up to now, there is no any algorithm or solver which can guarantee that the result obtained is the real global optimal solution, either in theory or in practice, so for PCC benchmark, it is possible that there might be better results in some cases that have not been found yet.

The five solvers, Baron, Antigone, Couenne, Lingo and Scip, are executed on the GAMS platform. In most cases, these solvers won't terminate until 1000 seconds or even more, even though the computation results have kept stable without any further improvements. Therefore, their actual computation times are taken, which is the time the best stable results reached firstly; the computation time for 1stOpt is also the time when the global result is first found; Matlab-GA takes their full process computation time.

The GOAs of the five solvers running on the GAMS platform are basically all based on the branch-and-bound deterministic algorithm, which has a distinctive feature that the final outcome remain stable and unchangeable. Under this circumstance, since the outcomes are not guaranteed to be the global optimal solution, once the result is unsatisfactory, there's no possibility of achieving better results if the configurations remain the same. Matlab-GA uses a genetic algorithm that contains lots of stochastic processes such as random selection, random mutation, etc., leading the computation results and time cost vary greatly and unsteadily with each run; therefore, the best result and the corresponding time among its ten runs are chosen and adopted. 1stOpt employs its own unique UGO algorithm, which is different from both deterministic and heuristic algorithms, although there are random processes in UGO algorithm, the calculation results remain stable in most cases, similar to Matlab-GA, the best result and corresponded time from its time runs are taken as the final solution.

The values of objective function, which is also referred as SSE defined in Eq. (1), from each solver in each problem are shown in Table 8. In addition to the values of SSE, the detailed parameter values corresponded are also given together for comparison purpose (see the appendix). In some cases, even though the objective function values yielded by different solvers differ slightly, the corresponding parameter values are completely different, which implies that these two outcomes are actually two totally different solutions. In this situation, if the objective function values are refereed only but ignoring the parameter values corresponded, it can result in a false judgment about whether the solution achieved is an approximate global optimal solution, or it is actually a fully different local solution.

The comparison results of the PCC benchmark test dataset are summarized below referring to the detailed of results in appendix:
1) Problem 1: the fitting function is a very simple power function with two parameters only, however, except 1stopt, all other solvers fail to get correct answer;
2) Problem 2: there are two global solutions, 1stOpt can get correct answer at 50% probability, while all others get different local solutions only;
3) Problem 3: 1stOpt gets the best solution with about 90% probability, while all others get different local solution only;
4) Problem 4: there are two global solutions, 1stOpt can gets correct answers for each time computations, while all others only find local solutions;
5) Problem 5: this is a very hard problem, 1stOpt gets the best solution in only 10% probability, all others can only obtain locally optimal solutions that differ very far;

6) Problem 6: 1stOpt gets the best solution with about 50% probability, while all others are local solutions;
7) Problem 7: 1stOpt gets the best solution with about 80% probability, while all others are local solutions;
8) Problem 8: the objective value of 2.71051E-20 form Antigone is already a pretty good solution and really tends to 0, however, compared to the value of 1.251867E-25 obtained by 1stOpt, since the parameter values of two solutions are completely different, strictly speaking, the result of Antigone is also local optimal solution; All others fail to find such solutions;
9) Problem 9: 1stOpt gets the best result with about 80% probability, the results of other solvers are all local optimization solutions.
10) Problem 10: Both Baron and 1stOpt get the best results, but the computation time of Baron is as much as 8 times that of 1stOpt; the results of other solvers are all local solutions.
11) Problem 11: Both Baron and 1stOpt have the best results and the same computation time; Antigone, although the objective value is near to the best one, is clearly not an approximate global optimal solution but a local one with referring to its parameter values; Couenne's results are obviously incorrect; the rest are all local solutions.
12) Problem 12: 1stOpt has the best result, the rests are all local solutions.
13) Problem 13: The results of Baron, Scip and 1stOpt are equally the best, but Baron and Scip are more efficient; the objective value of Couenne is only slightly different from the best solution, but according to its detailed parameter values, the solution is clearly not the global optimal solutions but a local one; the rests are local solutions.
14) Problem 14: The result of 1stOpt is the best, while others are local solutions.
15) Problem 15: Baron and 1stOpt have the same best results, and Baron is more efficient; the rests are local solutions.
16) Problem 16: 1stOpt gets the best solution, while Scip is the worst, and the rests are all local solutions.
17) Problem 17: 1stOpt gets the best solution, and the rests are all local solutions.
18) Problem 18: Baron and 1stOpt have the same best results, but Baron is more efficient, and the rest are all local solutions.
19) Problem 19: Baron and 1stOpt have the same best result, Baron is more efficient, and 1stOpt has only about 80% probability of successful, lower than that of Baron's 100%; the rest are all local solutions.
20) Problem 20: Baron and 1stOpt have the same best result, but 1stOpt has only 80% probability of successful, lower than that of Baron's 100%, the rests are all local solutions.
21) Problem 21: Baron and 1stOpt have the same best, 1stOpt is more efficient, but can only get that result in 80% success rate, lower than Baron's 100%; the rests are all local solutions.
22) Problem 22: 1stOpt has the best result, but it can only get this result in 50% success rate, and the remaining 50% probability yields a suboptimal result with objective function value of 186.141879, corresponding to parameter vector of b=[323.960345, 0.012968, 200001128069171.0, -0.239998070648, 843677740351.781, -0.237211], the suboptimal result is still better than the local solutions of all other solvers.
23) Problem 23: 1stOpt has the best result, but in only 90% success rate; the rests are all local

solutions.

24) Problem 24: 1stOpt has the best result, but in only 20% success rate; Antigone failed without giving any feasible solution, the rests are all local solutions.
25) Problem 25: Baron and 1stOpt have the same best result, Baron is more efficient in 100% success rate while 1stOpt can achieved global solution in only 50% success rate, the rests are all local solutions.
26) Problem 26: 1stOpt has the best result, but only in 90% success rate, and the other 10% probability yields a suboptimal local optimal solution which are same as the outcomes from Baron and Couenne, other results get the worse local solutions.
27) Problem 27: 1stOpt has the best result, but only in 40% probability, and the rests are all local solutions.
28) Problem 28: Both Couenne and 1stOpt can obtain the best results, 1stOpt is slightly more efficient, but it can only get this result in 40% probability, the rests are all local solutions.
29) Problem 29: Both Baron and 1stOpt have the same best results, 1stOpt is more efficient, but only in 70% probability of successful, lower than Baron's 100%; the rests are all local solutions.
30) Problem 30: 1stOpt has the best result, but can only be achieved with 50% probability; the objective values of Antigone and Couenne are the same and very close to the best one, but according to the corresponded parameter values, they should be local solutions; the rests are all local solutions.
31) Problem 31: Baron and 1stOpt have the same best result, Baron is more efficient and 1stOpt can only get the global solution with 50% probability, lower than Baron's 100%; the rests are all local solutions.
32) Problem 32: 1stOpt has the best result, but only with 40% probability, the other 60% probability all fall into the local solution of b=[-4.006080, -0.048869, -2.957393, -1.297652, -1.846793] with the objective value of 10.272014. This local solution is still better than the results obtained from all other solvers.
33) Problems 33: Baron and 1stOpt have the same best results, Baron is more efficient, 1stOpt can only get the result with 90% probability, lower than Baron's 100%; the rests are all local solutions.
34) Problem 34: 1stOpt has the best, but only with 50% probability; Baron and Scip have the same second-best result, which is also the local solution that 1stOpt falls into with another 50% probability; the rests are all local solutions.
35) Problem 35: 1stOpt is the best, but only with 60% probability; the rests are all local solutions.
36) Problem 36: Baron and 1stOpt are the same best, Baron is more efficient and 1stOpt can only get the global solution with 90% probability, lower than Baron's 100%; the rests are all local solutions.
37) Problem 37: Antigone, Couenne and 1stOpt have the same best results, 1stOpt is more efficient, but only with 70% probability, which is lower than Antigone and Couenne's 100%, and the rests are all local solutions.
38) Problem 38: 1stOpt is the best but in only 70% probability, the rests are all local solutions.

The comparison fitting results of each solver for each problem in PCC benchmark are shown in Fig. 5.

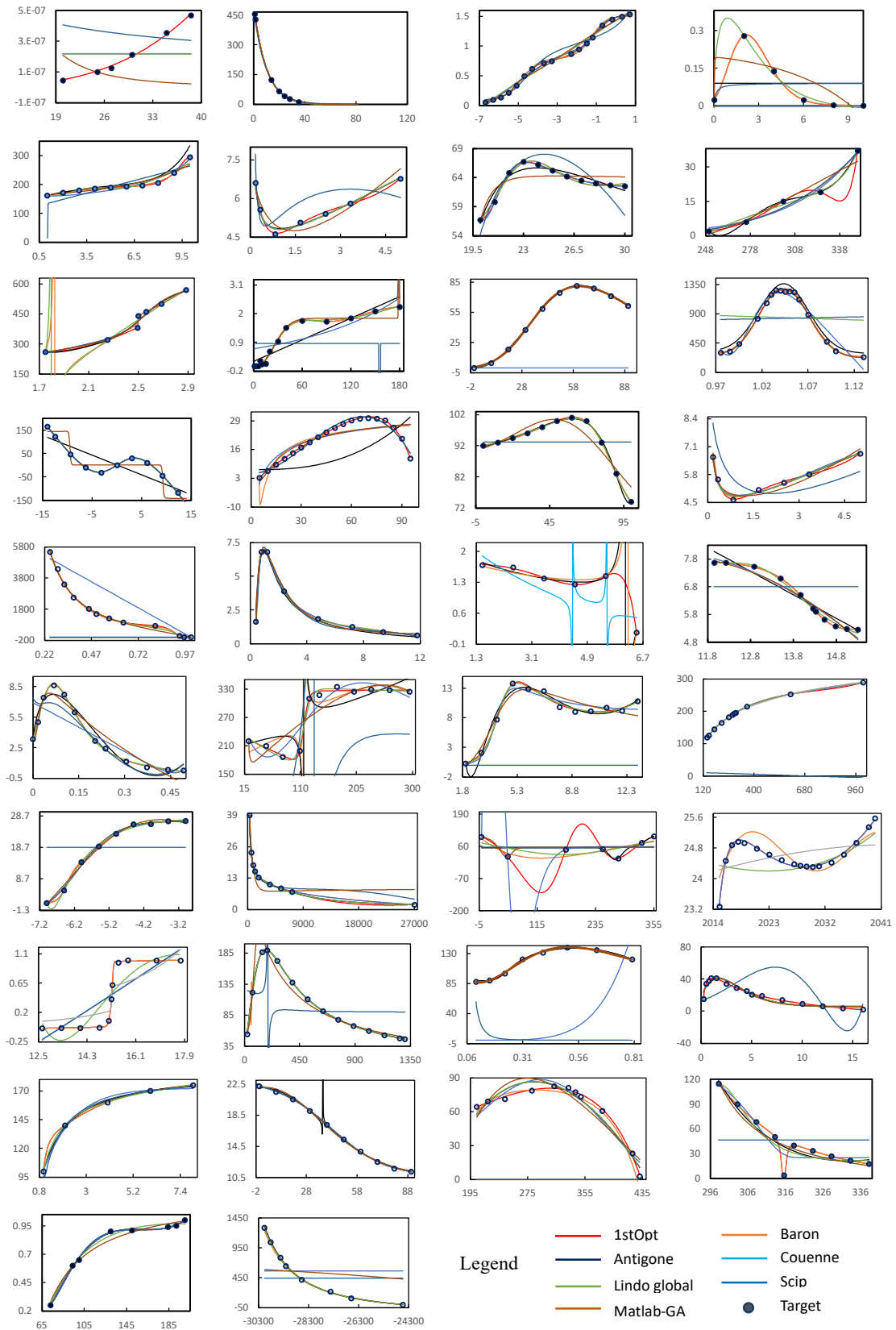

Fig. 5 Comparison results of each solver for PCC benchmark

Referring to the detailed calculation results in the appendix, the overall effect of each solver could be evaluated based on the success rate of obtaining the global solution of each solver, as shown in table 6, and refer to Fig. 6. Baron can get 13 global optimal solutions out of the 38 problems, with an overall success rate of 34.21%; Couenne has two problems solved with a success rate of 5.26%; both Antigone and Scip have only reached one global solution out of 38 problems with a success rate of 2.63%; Lingo and Matlab-GA failed at all, without even one problem solved successfully, especially for Matlab-GA, the results calculated each time are significantly different; 1stOpt performs the best in PCC benchmark, with 11 out of 38 problems being able to obtain the best optimal solution at 100% probability, and more than 62% successful probability for the remaining 27 problems, with an overall average success rate of 73.15%, which is more than twice higher than that of the second best performer in this study, Baron.

As for the efficiency, since it is only meaningful to compare the solvers that reach to the same global optimal solution in the same problem, considering that the success rate except for Baron and 1stOpt are extremely low in PCC benchmark, it is impossible to compare the efficiency of other solvers. Hence, the efficiency comparisons are made only for 13 cases where both Baron and 1stOpt obtain the best optimal solutions, see Table 7. It can be roughly estimated that overall average efficiency of 1stOpt is about twice higher than that of Baron, although for most individual problems, Baron is more efficient than 1stOpt.

**Table 6** Success Rate of Each Solver (%)

| Problem \ Solver | Baron | Antigone | Couenne | Lingo | Scip | Matlab-GA | 1stOpt |
|---|---|---|---|---|---|---|---|
| 1 | 0 | 0 | 0 | 0 | 0 | 0 | 100 |
| 2 | 0 | 0 | 0 | 0 | 0 | 0 | 50 |
| 3 | 0 | 0 | 0 | 0 | 0 | 0 | 90 |
| 4 | 0 | 0 | 0 | 0 | 0 | 0 | 100 |
| 5 | 0 | 0 | 0 | 0 | 0 | 0 | 10 |
| 6 | 0 | 0 | 0 | 0 | 0 | 0 | 50 |
| 7 | 0 | 0 | 0 | 0 | 0 | 0 | 80 |
| 8 | 0 | 0 | 0 | 0 | 0 | 0 | 80 |
| 9 | 0 | 0 | 0 | 0 | 0 | 0 | 80 |
| 10 | 100 | 0 | 0 | 0 | 0 | 0 | 100 |
| 11 | 100 | 0 | 0 | 0 | 0 | 0 | 100 |
| 12 | 0 | 0 | 0 | 0 | 0 | 0 | 100 |
| 13 | 100 | 0 | 0 | 0 | 100 | 0 | 100 |
| 14 | 0 | 0 | 0 | 0 | 0 | 0 | 100 |
| 15 | 100 | 0 | 0 | 0 | 0 | 0 | 100 |
| 16 | 0 | 0 | 0 | 0 | 0 | 0 | 100 |
| 17 | 0 | 0 | 0 | 0 | 0 | 0 | 100 |
| 18 | 100 | 0 | 0 | 0 | 0 | 0 | 100 |
| 19 | 100 | 0 | 0 | 0 | 0 | 0 | 80 |
| 20 | 100 | 0 | 0 | 0 | 0 | 0 | 80 |
| 21 | 100 | 0 | 0 | 0 | 0 | 0 | 80 |
| 22 | 0 | 0 | 0 | 0 | 0 | 0 | 50 |
| 23 | 0 | 0 | 0 | 0 | 0 | 0 | 90 |
| 24 | 0 | 0 | 0 | 0 | 0 | 0 | 20 |
| 25 | 100 | 0 | 0 | 0 | 0 | 0 | 50 |
| 26 | 0 | 0 | 0 | 0 | 0 | 0 | 90 |
| 27 | 0 | 0 | 0 | 0 | 0 | 0 | 40 |
| 28 | 0 | 0 | 100 | 0 | 0 | 0 | 40 |
| 29 | 100 | 0 | 0 | 0 | 0 | 0 | 70 |
| 30 | 0 | 0 | 0 | 0 | 0 | 0 | 50 |
| 31 | 100 | 0 | 0 | 0 | 0 | 0 | 50 |
| 32 | 0 | 0 | 0 | 0 | 0 | 0 | 40 |
| 33 | 100 | 0 | 0 | 0 | 0 | 0 | 90 |
| 34 | 0 | 0 | 0 | 0 | 0 | 0 | 50 |
| 35 | 0 | 0 | 0 | 0 | 0 | 0 | 60 |
| 36 | 100 | 0 | 0 | 0 | 0 | 0 | 80 |
| 37 | 0 | 100 | 100 | 0 | 0 | 0 | 70 |
| 38 | 0 | 0 | 0 | 0 | 0 | 0 | 60 |
| Average Success Rate (%) | 34.21 | 2.63 | 5.26 | 0.00 | 2.63 | 0.00 | 73.15 |

**Table 7** Average Time Used of Baron and 1stOpt for Successful Problems

| Problem | | 10 | 11 | 13 | 15 | 18 | 19 | 20 | 21 | 25 | 29 | 31 | 33 | 36 | total |
|---|---|---|---|---|---|---|---|---|---|---|---|---|---|---|---|
| Time (s) | Baron | 34.0 | 15.0 | 1.0 | 1.0 | 1.5 | 1.5 | 1.0 | 78.8 | 1.0 | 6.0 | 2.0 | 4.0 | 1.0 | 147.8 |
| | 1stOpt | 4.0 | 15.0 | 2.5 | 3.0 | 4.0 | 2.5 | 2.0 | 5.0 | 8.0 | 2.0 | 15.0 | 10.0 | 2.0 | 75 |

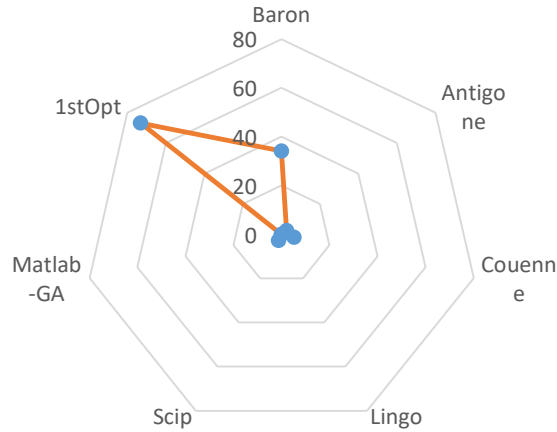

Fig. 6 Success Rate of each Solver（%）

**Table 8** Objective Function Values of PCC Benchmark Problems from Each Solver

| Solver / Problem | Baron | Antigone | Couenne | Lingo | Scip | Matlab-GA | 1stOpt |
|---|---|---|---|---|---|---|---|
| 1 | 1.333068E-13 | 1.33307E-13 | 1.33307E-13 | 1.33307E-13 | 2.97525E-13 | 3.593366E-13 | 5.938663E-16 |
| 2 | 0.00299032 | 0.00391558 | 0.00688819 | 0.00431527 | 4.665052E+2 | 0.007932975 | 0.000745813 |
| 3 | 0.01005775 | 0.01590848 | 0.01590848 | 0.04949805 | 0.11571187 | 0.04149313 | 0.000916735 |
| 4 | 0.00000011 | 0.05778266 | 0.09833024 | 0.00061758 | 0.05938093 | 0.025760299 | 1.5306075E-26 |
| 5 | 60.74644712 | 25.55889645 | 1.015841E+3 | 9.057286E+2 | 2.518579E+4 | 1314.3774608 | 2.751590632 |
| 6 | 0.06097893 | 0.07899019 | 0.05843872 | 0.06097754 | 3.84683560 | 0.611337649 | 0.007128369 |
| 7 | 1.54802194 | 7.50860528 | 0.07730189 | 1.39222663 | 57.9845223 | 24.5993411 | 0.00389328 |
| 8 | 15.16351237 | 2.71051E-20 | 16.74389478 | 6.46699570 | 17.40494834 | 46.33834746 | 1.251867E-25 |
| 9 | 1191.5286 | 1.723727E+3 | 1.864564E+3 | 1971.28317 | 1.864564E+3 | 1.77714e+03 | 0.3666920 |
| 10 | 0.06654203 | 1.82332649 | 4.08569134 | 0.07865553 | 10.97874268 | 0.16307029 | 0.06654202 |
| 11 | 0.00704569 | 0.00729984 | 3.258914E+4 | 0.01864764 | 0.01864764 | 4.43442032 | 0.007045685 |
| 12 | 5.634793E+3 | 7.509296E+3 | 1.440851E+5 | 2.571216E+6 | 2.586745E+6 | 6142.972817 | 1391.501878 |
| 13 | 24.83018843 | 2.117121E+4 | 24.99678686 | 33.35848010 | 24.83018364 | 4190.433475 | 24.83018363 |
| 14 | 4.067720E+2 | 33.82722561 | 4.508885E+2 | 25.25308277 | 25.25308388 | 411.5534198 | 0.909227034 |
| 15 | 0.29936497 | 1.75498593 | 9.607825E+4 | 0.98870723 | 6.586250E+2 | 68.55265925 | 0.299364969 |
| 16 | 0.03495035 | 0.06097539 | 0.09270095 | 0.06097754 | 7.07618234 | 0.28979101 | 0.00712837 |
| 17 | 4.301794E+4 | 7.555325E+7 | 1.694418E+7 | 7.974938E+4 | 7.469755E+7 | 158654.0395 | 34105.93604 |
| 18 | 0.00551218 | 0.00693548 | 0.00885065 | 0.01876434 | 0.13393480 | 0.11187831 | 0.00551218 |
| 19 | 0.0012525071 | 0.01661039 | 0.26418161 | 0.012957474 | 0.26418161 | 0.02702564 | 0.0012525071 |
| 20 | 0.00086491 | 0.56099429 | 0.35950147 | 0.05563105 | 9.56709091 | 0.41550289 | 0.00086490 |
| 21 | 0.1814943 | 4.1689152 | 33.624865 | 0.2106728 | 26.9151058 | 7.7368775 | 0.1814943 |
| 22 | 3.399997E+3 | 5.885228E+3 | 5.567073E+3 | 483.32350 | 8.108606E+5 | 8378.444731 | 179.3335177 |
| 23 | 2.12289246 | 5.68622712 | 8.14877184 | 2.12333274 | 1.132110E+3 | 17.86916266 | 1.978953145 |
| 24 | 0.99685339 | Failed | 0.99685339 | 3.6361175 | 3.936153E+5 | 40.36491692 | 0.52007113 |
| 25 | 0.44499232 | 2.79312763 | 2.79312763 | 1.43570424 | 7.891489E+2 | 6.80895741 | 0.444992322 |
| 26 | 0.08897309 | 0.62251302 | 0.08897309 | 0.62251342 | 1.35478951 | 78.40704417 | 0.005260412 |
| 27 | 3.289399E+3 | 3.291632E+3 | 4.317044E+3 | 4.375858E+3 | 7.646795E+3 | 4433.438896 | 6.681978318 |
| 28 | 0.65965216 | No solution | 0.02493958 | 2.84091799 | 1.090121E+4 | 3.65573320 | 0.02493958 |
| 29 | 0.00073779 | 0.70660012 | 0.72451259 | 0.49538789 | 0.72451259 | 0.38729520 | 0.000737787 |
| 30 | 0.66930130 | 0.46183999 | 0.46183999 | 0.651761430 | 2.476439E+4 | 999.2931747 | 0.461594701 |

| 31 | 2.24560393 | 5.38894747 | 8.214714E+4 | 13.93232765 | 1.000686E+5 | 17.07814305 | 2.24560393 |
| 32 | 1.742015E+2 | 1.742015E+2 | 1.742015E+2 | 1.742015E+2 | 7.084184E+3 | 183.6248605 | 8.31641750 |
| 33 | 2.46322E-20 | 6.89358250 | 39.02123764 | 0.84093052 | 1.090936E+5 | 8.64425832 | 5.16987E-26 |
| 34 | 0.01228654 | 0.14654921 | 0.01237477 | 0.01237477 | 0.01228654 | 0.35758658 | 0.010808142 |
| 35 | 75.57370284 | 3.296046E+2 | 4.625350E+2 | 329.78055817 | 6.570858E+3 | 814.3733233 | 30.02585461 |
| 36 | 0.81336718 | 1.461242E+3 | 1.178249E+3 | 1330.351267 | 1.085889E+4 | 1606.349399 | 0.81336718 |
| 37 | 0.00085308 | 0.00069673 | 0.00069673 | 0.00388164 | 0.00085308 | 0.00868975 | 0.00069673 |
| 38 | 7.176911E+2 | 7.176883E+2 | 1.475831E+6 | 7.165598E+2 | 1.584736E+6 | 1199841.698 | 676.71429733 |
| No. of Success | 13 | 1 | 2 | 0 | 1 | 0 | 38 |

## 5 Discussion and Conclusion

A new challenging curve fitting test dataset of PCC benchmark has been proposed. Based on NIST and PCC benchmark, seven of today's lead global optimization solvers were checked and compared in terms of effect and efficiency in the area of nonlinear curve fitting. All test problems are characterized by small-scale with the maximum parameter numbers no more than 8, and the larger-scale problems are not included and discussed in this study. When evaluating GOAs or relative solvers, while success with small-scale optimization problems doesn't necessarily ensure success with large-scale ones, it is certain that a GOA or solver performs poorly on small-scale problems, it won't fare well with large-scale problems. Therefore, the problems with small-scale and free search domain but high difficulty are more suitable for testing and verifying the effectiveness of GOAs or relative solvers. In addition, while all the solvers demonstrate an overall good performance in NIST benchmark, their comparative results in PCC benchmark differ significantly from each other and often fail to find global solution most of the time, except for 1stOpt, which is able to find the certified solution in NIST benchmark every time and obtain the best solution for most problems in PCC benchmark. This suggests that the PCC benchmark is more distinguishable, challenging, and more suitable in testing the performance of GOAs in the field of nonlinear curve fitting.

There are three points that need special attention. The first one is that the proposed PCC benchmark test dataset are very complicated global optimization problems, it is currently impossible to theoretically prove that the best results obtained up to now in this study are the true global optimization solutions, and there might exist better results that are found in the future. Secondly, each solver comes along with various configuration options, different choices of these settings may lead to different results, even the same solver in different versions may produce varying outcomes. Taking 1stOpt as an example, if adding one statement "Hardness = 2" (the default value is 1), the average success percentage rate will increase about 20%, in the meanwhile, consuming more time. In this study, the default settings are used for all GOA solvers. Thirdly, since this study primarily focuses on nonlinear curve fitting, the results and conclusions obtained may not necessarily apply to other types of global optimization problems.

The following preliminary conclusions can be concluded:
1) The PCC benchmark test dataset proposed in this paper is in low dimension but high difficulty, and is challenging enough to test, verify, and distinguish the performance of current GOAs. Furthermore, it can also be a reliable resource for developers to improve their algorithms.
2) In terms of effect, for the NIST benchmark, 1stOpt, Baron, Lingo, Antigone, Couenne,

Scip and Matlab-GA are ranked from best to worst; For the PCC benchmark, 1stOpt performs the best, Baron and Couenne are second and third, respectively. Antigone and Scip tie for fourth, while Lingo and Matlab-GA perform the worst;

3) In terms of efficiency, for the NIST benchmark, the order from high-efficient to low-efficient is Baron, Couenne, Scip, Antigone, 1stOpt, Matlab-GA and Lingo; For the PCC benchmark, Antigone, Couenne, Scip, Lingo and Matlab-GA cannot be counted because of their extremely low success rates, and for the remained 1stOpt and Baron, the former is about twice higher than the latter.

4) The NIST benchmark is mainly designed to test LOA or local optimization solvers, as it is relatively too simple and thus it's not well-suited for testing GOA or global optimization solvers.

5) During the study process, it is not enough to only consider the objective function values, but also need to take corresponding parameter values into account; otherwise, it will be impossible to exactly judge whether the obtained result is the true global optimal solution or a local one.

The PCC benchmark is a unique test dataset and can be employed for efficiently test the performance of a GOA or relative solver within a short time. In the future, it is necessary to test GOA or relevant solver with the more large-scale test problems, and expand from the field of nonlinear curve fitting to other areas of global optimization.

**Data Availability** All the data are reported in the paper.

**Funding** The authors declare that no funds, grants, or other support were received during the preparation of this manuscript.

**Conflict of Interests** The authors have no relevant financial or non-financial interests to disclose.

# References

1. Belotti, P., Lee, J., Liberti, L., Margot, F., Wächter, A.: Branching and bounds tightening techniques for non-convex MINLP. Optim. Method. Softw. 24(4–5), 597–634 (2009)
2. Bestuzheva, K., Besançon, M., Chen, W.K., Chmiela, A., Donkiewicz, T., van Doornmalen, J., Eifler, L., Gaul, O., Gamrath, G., Gleixner, A., et al., 2021. The SCIP optimization suite 8.0. arXiv preprint arXiv:2112.08872.
3. Dorigo, M., Birattari, M., & Stutzle, T. (2006). Ant colony optimization. IEEE computational intelligence magazine, 1(4), 28–39.
4. E. Goodarzi, M. Ziaei, and E. Z. Hosseinipour, "Optimization analysis using LINGO and MATLAB," in Introduction to Optimization Analysis in Hydrosystem Engineering, vol. 25, pp. 149–193, Springer, 2014.
5. F. Glover, Future paths for integer programming and links to artificial intelligence, Comp. Oper. Res.13(1986)533-549.
6. Gergel, V., Sergeyev, Y.: Sequential and parallel algorithms for global minimizing functions with lipschitzian derivatives. Computers and Mathematics with Applications 37(4), 163–179 (1999). DOI https://doi.org/10.1016/S0898-1221(99)00067-X. URL https://www.sciencedirect.com/science/article/pii/S089812219900067X
7. Ghildyal V., Sahinidis N.V. (2001) Solving Global Optimization Problems with Baron. In: Migdalas A., Pardalos P.M., Värbrand P. (eds) From Local to Global Optimization. Nonconvex Optimization and Its Applications, vol 53. Springer, Boston, MA. https://doi.org/10.1007/978-1-4757-5284-7_10
8. Hey, A.M.: Towards Global Optimisation 2. Journal of the Operational Research Society 30(9), 844 (1979). DOI 10.1057/jors.1979.201. URL https://doi.org/10. 1057/jors.1979.201
9. J. Kennedy and R. Eberhart, "Particle swarm optimization," Proceedings of ICNN'95 - International Conference on Neural Networks, Perth, WA, Australia, 1995, pp. 1942-1948 vol.4, doi: 10.1109/ICNN.1995.488968.
10. J.H. Holland (1975) Adaptation in Natural and Artificial Systems, University of Michigan Press, Ann Arbor, Michigan; re-issued by MIT Press (1992).
11. v56. Jones, D.R., Perttunen, C.D., Stuckman, B.E.: Lipschitzian optimization without the Lipschitz constant. Journal of Optimization Theory and Application 79(1), 157–181 (1993). DOI 10.1007/BF00941892
12. Khajavirad, A. and N. V. Sahinidis, A hybrid LP/NLP paradigm for global optimization relaxations, Mathematical Programming Computation, 10, 383-421, 2018.
13. Kirkpatrick, S., Gelatt, Jr., C.D., Vecchi, M.P.: Optimization by simulated annealing. Science, 220, 671–680 (1983)
14. Liu, J., Ploskas, N. & Sahinidis, N.V. Tuning BARON using derivative-free optimization algorithms. J Glob Optim 74, 611–637 (2019). https://doi.org/10.1007/s10898-018-0640-3
15. MathWork, Inc. Global Optimization Toolbox, https://www.mathworks.com/help/gads/
16. M. Dorigo and T. Stützle, "Ant Colony Optimization: Overview and Recent Advances," in Handbook of Metaheuristics, vol. 272, M. Gendreau, J.-Y. Potvin, Eds., Springer, Cham, 2019, pp. 311-351. [Online]. Available: https://doi.org/10.1007/978-3-319-91086-4_10
17. Misener, R., Floudas, C.A. ANTIGONE: Algorithms for coNTinuous / Integer Global Optimization of Nonlinear Equations. J Glob Optim 59, 503–526 (2014).

https://doi.org/10.1007/s10898-014-0166-2

18. Misener, R., Floudas, C.A. GloMIQO: Global mixed-integer quadratic optimizer. J Glob Optim 57, 3–50 (2013). https://doi.org/10.1007/s10898-012-9874-7
19. N. Sahinidis and M. Tawarmalani, "GLOBAL OPTIMIZATION WITH GAMS/BARON", 2014.
20. Nikolaos V. Sahinidis (2009) Global optimization, Optimization Methods and Software, 24:4-5, 479-482, DOI: 10.1080/10556780903135287
21. NIST: StRD Nonlinear Least Squares Regression Datasets. https://www.itl.nist.gov/div898/strd/nls/nls_main.shtml
22. Poli, R., Kennedy, J. & Blackwell, T. Particle swarm optimization - An overview. Swarm Intell 1, 33–57 (2007). https://doi.org/10.1007/s11721-007-0002-0.
23. Shubert, B.O.: A sequential method seeking the global maximum of a function. SIAM Journal on Numerical Analysis 9, 379–388 (1972). DOI 10.1137/0709036
24. Stefan Vigerske & Ambros Gleixner (2018) SCIP: global optimization of mixed-integer nonlinear programs in a branch-and-cut framework, Optimization Methods and Software, 33:3, 563-593, DOI: 10.1080/10556788.2017.1335312
25. Stripinis, L., Paulavičius, R.: DIRECTGOLib - DIRECT Global Optimization test problems Library (2022). DOI 10.5281/zenodo.6617799. URL https://doi.org/10. 5281/zenodo.6617799
26. Xianyun Cheng, Fuxin Chai, Jing Gao, Kelei Zhang. 1stOpt and Global Optimization Platform——Comparison and Case Study. Proceedings of 2011 4th IEEE International Conference on Computer Science and Information Technology (ICCSIT 2011) VOL. 04

# Appendix: Detailed Results of PCC Benchmark Test Dataset

The appendix gives the detailed results from each solver of all 36 problems in the PCC Benchmark test dataset, including the objective function values, parameter values and time used, the rows with dark background colors in the table represent the best solution obtained by the corresponding solvers.

Appx. Table 1: Detailed Results of Problem 1

| Solver | Obj. Values | 2 Parameters | Time (s) |
|---|---|---|---|
| Baron | 1.333068E-13 | 0.21903333E-006, 0.9735653921E-013 | 0.1 |
| Antigone | 1.33307E-13 | 0.00000022, 0.00000000 | 0.1 |
| Couenne | 1.33307E-13 | 0.00000022, 0.00000002 | 0.1 |
| Lingo | 1.33307E-13 | 0.00000022, 0.00000000 | 0.1 |
| Scip | 2.97525E-13 | 0.00000150, -0.43631118 | 0.2 |
| Matlab-GA | 3.593366E-13 | 0.00522131, -3.37873144 | 1.5 |
| 1stOpt | 5.938663E-16 | 1.59002014E-12, 3.45494899 | 0.1 |

Appx. Table 2: Detailed Results of Problem 2

| Solver | Obj. Values | 6 Parameters | Time (s) |
|---|---|---|---|
| Baron | 0.00299032 | -1.60803, -4.99863E+2, -0.100267, 4.93408E+3, -1.422547, 1.8405E+2 | 3.0 |
| Antigone | 0.00391558 | -1.61605, -4.99906E+2, -0.100115, -5.28469E+6, -3.38096, 1.5719E+3 | 10.0 |
| Couenne | 0.00688819 | -1.623335, -0.712666, -0.100107, 0.886066, -0.00000869 -1.73346776 | 4.3 |
| Lingo | 0.00431527 | -1.613085, -5.588210, -0.1013398, -0.00008276, 0.9129759, -0.777444 | 2.0 |
| Scip | 4.665052E+2 | 0.913187, -2.47673E+1, -0.094116, 0.1486097, 0.02017663, -0.811161 | 4.5 |
| Matlab-GA | 0.007932975 | -1.5927888, -499.20744, -0.1003532, 689.8065, -0.446666, -32.457549 | 16.0 |
| 1stOpt | 0.000745813 (50% Pr.) | [-1.519158, -542.37158, -0.099986, 0.0002517, 2.444117, -13.210135], [-1.519158,-499.7221, -0.099986, -592291.7768, -2.444118, 11.21013] | 7.0 |

Appx. Table 3: Detailed Results of Problem 3

| Solver | Obj. Values | 6 Parameters | Time (s) |
|---|---|---|---|
| Baron | 0.01005775 | 0.7326377, 0.0188533, 3.5886105, 0.3834541, 0.00000139, 3.064187 | 1.0 |
| Antigone | 0.01590848 | -0.4843523, 9.610936, -1.5463586, 0.7858020, 0.0366302, 0.6767113 | 122 |
| Couenne | 0.01590848 | -0.4843523, 9.610942, -1.54635896, 0.785802, 0.0366302, 0.6767113 | 1.5 |
| Lingo | 0.04949805 | -0.75358556, -1.520213, 0.049488, -0.00757813, 0.00003073, 64.7367 | 0.5 |
| Scip | 0.11571187 | -0.778474, 7.089495, 0.81466366, 0.2876983, -0.9579807, -0.033699 | 1.0 |
| Matlab-GA | 0.04149313 | -0.16432536, 477.0499, 3.9550722, 0.987716, 0.3888601, 0.45152857 | 11.0 |
| 1stOpt | 0.000916735 (90% Pr.) | 0.698457, 6.143456E-10, 6.3195982, 0.420541, 0.1720433, 2.433167 | 6.0 |

Appx. Table 4: Detailed Results of Problem 4

| Solver | Obj. Values | 6 Parameters | Time (s) |
|---|---|---|---|
| Baron | 0.00000011 | 0.0011895, 0.119137, 1.471494, 7.799656, 7.400285E+4, 4.442459E+7 | 3.7 |
| Antigone | 0.05778266 | 1.0000E+6, 0.000000, 44.027819, -4.6334846, -8.9400E+3, -1.000E+5 | 5.0 |

| Couenne | 0.09833024 | 0.000000, 0.000000, 1.77636E-15, 0.000000, 0.000000, -0.49999999 | 1.0 |
|---|---|---|---|
| Lingo | 0.00061758 | 0.630024, 0.578827, 0.00147435, 0.5000, -1.44463E+8, 8.807665E+9 | 18.0 |
| Scip | 0.05938093 | 0.02321082, 0.00000, 0.89233118, -0.366352, 0.7366682, 8.186676 | 1.0 |
| Matlab-GA | 0.025760299 | -369.034254, 1.2718887, 63.789170, -2.2343045, 94.21665, 409.132 | 10.0 |
| 1stOpt | 1.5306075E-26 | [0.0022116, -0.1282624, 1.385715, 7.33447, 8.57023E-5, -0.0536963], [0.0022116, 0.1282624, 1.385715, 7.33447, 8.57023E-5, -0.0536963] | 20 |

Appx. Table 5: Detailed Results of Problem 5

| Solver | Obj. Values | 6 Parameters | Time (s) |
|---|---|---|---|
| Baron | 60.74644712 | 2.30748E+8, -3.0972E+3, -0.461796, 1.37617E+6, -0.170895, -8.0639 | 1.5 |
| Antigone | 25.55889645 | 0.0011582, -0.99999, -0.0003777, 0.0003752, 1.002964, -4.59367E+1 | 2.0 |
| Couenne | 1.015841E+3 | -1.400E+1, 0.0098947, -0.1732868, -1.10083278, 0.003232, 2.268599 | 16.0 |
| Lingo | 9.057286E+2 | -1.1343E+12, -5.59485E+9, -0.518279, 1.04261E+9, 2.43023, 27.1325 | 3.0 |
| Scip | 2.518579E+4 | 1.19207E+2, 0.00000, 1.30343E+4, 8.8713E+6, -3.5253E+2, 14.4588 | 1.0 |
| Matlab-GA | 1314.3774608 | 146.513295, 150.42967, 9.8764265, -0.1735306, 0.567057, -13.95439 | 4.2 |
| 1stOpt | 2.751590632 (10% Pr.) | 153.2675, 0.00059321, -0.9574511, -4.799846E-9, 9.262531, 9.039135 | 30 |

Appx. Table 6: Detailed Results of Problem 6

| Solver | Obj. Values | 6 Parameters | Time (s) |
|---|---|---|---|
| Baron | 0.06097893 | 3.378793, -2.3512E+2, -0.0026883, 2.35327E+2, -0.999771, 0.637151 | 0.5 |
| Antigone | 0.07899019 | 2.32539E+8, 0.051086, -0.0173788, -1.7903E+1, 0.0067516, 0.656033 | 1.5 |
| Couenne | 0.05843872 | 3.433563, 0.435675, -0.00194228, -0.24423888, -1.163909, 0.6347529 | 150 |
| Lingo | 0.06097754 | 3.378777, -3.79688E+2, -0.002679, 3.798948E+2, -0.999858, 0.637191 | 1.0 |
| Scip | 3.84683560 | 1.14107E+3, 0.333142, -0.0167729, -5.245502, -0.1537995, -3.208953 | 3.0 |
| Matlab-GA | 0.611337649 | 22.95958, -0.0248411, 0.1698506, -1.6531279, 0.2104147, 1.2402976 | 32 |
| 1stOpt | 0.007128369 (50% Pr.) | -1.103458E-7, -0.062913, 0.0103969, 21.45572, 0.0451212, 217.97633 | 31 |

Appx. Table 7: Detailed Results of Problem 7

| Solver | Obj. Values | 6 Parameters | Time (s) |
|---|---|---|---|
| Baron | 1.54802194 | 39.27304, 1.00E+10, -2.52723E+3, -0.296194, -0.782884, 0.7765362 | 1.5 |
| Antigone | 7.50860528 | 88.67655, -1.46758E+3, 6.703762, -1.0269E-11, -0.565588, -0.895752 | 460 |
| Couenne | 0.07730189 | 1.28519E+2, 2.995653, -3.7568E+11, -1.270304, 0.0885004, -3.618028 | 30.0 |
| Lingo | 1.39222663 | 40.274686, 1.250E+9, -5.34681E+3, -0.338129, -0.7166243, 0.7418804 | 29.0 |
| Scip | 57.9845223 | 1.642108E+2, -1.52716E+1, 34.0070, -0.189542, 0.052501, -0.797417 | 7.0 |
| Matlab-GA | 24.5993411 | 65.49948, -66.11828, 168.2011, -0.17564897, -0.3547658, -0.0463128 | 4.5 |
| 1stOpt | 0.00389328 (80% Pr.) | 48.6653, 3385.108, -7438273877.457, -1.052315, -0.248421, 0.396466 | 4.0 |

Appx. Table 8: Detailed Results of Problem 8

| Solver | Obj. Values | 5 Parameters | Time (s) |
|---|---|---|---|
| Baron | 15.16351237 | 0.06368604, 0.01841413, 2.580498E+4, -2.47131E+5, -3.75806394 | 4.5 |
| Antigone | 2.71051E-20 | 8.817284E+5, -0.00340865, 8.819610E+5, -2.92314E+2, -6.49983E+5 | 76.0 |
| Couenne | 16.74389478 | 0.01465789, 0.02241918, -0.00012564, -4.48746E+3, -0.61280977 | 54.0 |

| Solver | Obj. Values | 5 Parameters | Time (s) |
|---|---|---|---|
| Lingo | 6.46699570 | 1.25885E-14, 0.09915911, -5.03004E+3, 3.58516943, 5.103937E+3 | 1.5 |
| Scip | 17.40494834 | 0.61573812, -0.20406817, 2.580498E+4, -2.30074E+3, 0.11264082 | 1.0 |
| Matlab-GA | 46.33834746 | 7776.13393, -0.021297289, -286.317898, 128.05021808, 440.754761 | 36 |
| 1stOpt | 1.251867E-25 (80% Pr.) | 1.11507877E-7, 0.065751, -4246745414541.75, -7737.47846, 0.618599 | 50.0 |

Appx. Table 9: Detailed Results of Problem 9

| Solver | Obj. Values | 6 Parameters | Time (s) |
|---|---|---|---|
| Baron | 1191.52859 | 142.04000, -6.25410, -40.87474, -89098.74499, -355.08261, 334.37932 | 0.5 |
| Antigone | 1.723727E+3 | 11.06580, 3.32621, -1.3208E+2, -0.009648, -2.79855E+1, 2.065587E+2 | 1.0 |
| Couenne | 1.864564E+3 | -1.04526E+1, 12.88994, 1.69229E+2, -0.993329, 26.64630, 1.29674E+2 | 11.0 |
| Lingo | 1971.28317 | -1.114108, -0.155708, -0.90966, -1.322834, -548.57442, 391.483335 | 1.0 |
| Scip | 1.864564E+3 | -1.04526E+1, 12.88994, 1.69229E+2, -0.993329, 26.64630, 1.29674E+2 | 1.0 |
| Matlab-GA | 1.77714e+03 | -9.648325, 2.33937, 1.85718e+02, 0.18096, 24.530841, 1.250866e+02 | 9.0 |
| 1stOpt | 0.36669200 (80% Pr.) | 1365.39674, 0.008109, -58.79124, -1.01059, -3409.99558, 789.31695 | 5.0 |

Appx. Table 10: Detailed Results of Problem 10

| Solver | Obj. Values | 7 Parameters | Time (s) |
|---|---|---|---|
| Baron | 0.06654203 | 2.33075928, -2.25324443, -0.00095946, 3.07215727, -0.63459110, -0.00028578, 91.23020201 | 34.0 |
| Antigone | 1.82332649 | -1.33733E+4, 5.716890E+3, 3.56488E-12, -1.50969E+5, 6.668214E+3, 3.35341E-12, -1.47471E+5 | 1.0 |
| Couenne | 4.08569134 | 0.00000000, 0.76163414, -2.00000000, 1.785858E+2, 0.00002182, 0.00000133, -2.78400E+3 | 1.0 |
| Lingo | 0.07865553 | 6.203525E+3, -6.20178E+3, -8.17912E-9, 75.82557436, -2.00401632, -0.00117016, 2.91544137 | 1.0 |
| Scip | 10.97874268 | 0.86313467, -4.23348E+1, -5.79610047, 1.555727E+2, -6.52202E+3, -6.00000E+2, -8.00000E+2 | 1.0 |
| Matlab-GA | 0.16307029 | 1.81939856, 34.43146047, -6.56578017, 179.18632377, -1.80416000, -0.00134009, 5.405509675486806 | 2.0 |
| 1stOpt | 0.06654202 | 2.33075927, -2.25324442, -0.00095946, 3.07215726, -0.63459109, -0.00028577, 91.23020201, | 4.0 |

Appx. Table 11: Detailed Results of Problem 11

| Solver | Obj. Values | 7 Parameters | Time (s) |
|---|---|---|---|
| Baron | 0.00704569 | -2.74112E+2, 25.45637012, 2.741512E+2, -0.09894504, 19.23286049, 0.01740064, -0.67053062 | 15.0 |
| Antigone | 0.00729984 | 2.625293E+6, -2.49777E+5, -2.62578E+6, -0.10198578, -2.01777E+5, 0.01705535, 6.825573E+3 | 69.0 |
| Couenne | 3.258914E+4 | 0.00000000, 0.00000000, 0.00000000, 0.00000000, 0.00000000, 0.00000000, 0.00000000 | 1.0 |
| Lingo | 0.01864764 | -4.27401E+1, -2.90239558, 42.74192082, 0.06352225, 0.11210860, 0.08316135, 0.06016641 | 1.0 |
| Scip | 0.01864764 | -4.27401E+1, -2.90239558, 42.74192082, 0.06352225, 0.11210860, 0.08316135, 0.06016641 | 1.0 |
| Matlab-GA | 4.43442032 | 909.398702, 96.8342394, -738.058802, 0.02275045, 321.25258464, -0.06821197, 0.68506228 | 6.0 |
| 1stOpt | 0.007045685 | -274.112079, 25.45637116, 274.151612, -0.09894482, -0.67053013, 19.23284579, 0.017400647 | 15.0 |

Appx. Table 12: Detailed Results of Problem 12

| Solver | Obj. Values | 8 Parameters | Time (s) |
| --- | --- | --- | --- |
| Baron | 5.634793E+3 | 2.113840E+3, -1.56799E+3, -4.65043E+5, 18.59414811, 1.04008545, 4.658367E+5, 18.61606476,1.04010431 | 60.0 |
| Antigone | 7.509296E+3 | 1.792470E+5, -1.47605E+2, 1.057384E+3, 26.26556134, 1.04377024, -1.79667E+5, 0.00000007, -1.00000E+6 | 2.0 |
| Couenne | 1.440851E+5 | 9.209378E+3, -1.51625E+4, 7.973932E+3, 0.00000047, -4.31830E+1, -2.05465E+3, -2.45373E+1, 010.97711556 | 10.0 |
| Lingo | 2.571216E+6 | 4.390761E+2, -4.56087E+2, 4.390761E+2, 0.00000000, 0.00000000, 4.390761E+2, 0.00000000, 0.00000000 | 1.0 |
| Scip | 2.586745E+6 | 3.082181E+2, 2.717571E+2, 1.193267E+2, 0.00000000, 0.00000000, 1.193267E+2, 0.00000000, 0.00000000 | 1.0 |
| Matlab-GA | 6142.972817 | -1195.251006, 1248.330659, 503.5691175, -22.26412009, 0.93263219, 1180.1472491, 24.40550212, 1.04253165 | 25.0 |
| 1stOpt | 1391.501878 | 466.131697, -205.728688, 1082.883499, -27.64683369, 1.04396274, -95.64724535, -148.59051485, 1.04468788 | 8.5 |

Appx. Table 13: Detailed Results of Problem 13

| Solver | Obj. Values | 5 Parameters | Time (s) |
| --- | --- | --- | --- |
| Baron | 24.83018843 | -5.76532E+3, 1.155179E+4, 0.03910648, 0.00998102, -0.36340820 | 1.0 |
| Antigone | 2.117121E+4 | 1.71035827, -1.69915681, -5.62500E+5, -0.00001761, 0.00001761 | 1.0 |
| Couenne | 24.99678686 | -7.58272E+3, 1.514634E+4, 1.02975352, 0.36379981, 0.00708805 | 14.0 |
| Lingo | 33.35848010 | -7.32422E+6, 1.464846E+7, 0.00003183, 0.00000771, -0.35529737 | 660.0 |
| Scip | 24.83018364 | -5.78162E+3, 1.158439E+4, 0.96100471, -0.36341288, 0.00995152 | 1.0 |
| Matlab-GA | 4190.433475 | 0.93022734, 1.90764562, 75.26010117, -9.44321620, 9.29816318 | 2.0 |
| 1stOpt | 24.83018363 | -5781.625354, 11584.4006729, 0.03899525, 0.00995151, -0.36341288 | 2.5 |

Appx. Table 14: Detailed Results of Problem 14

| Solver | Obj. Values | 5 Parameters | Time (s) |
| --- | --- | --- | --- |
| Baron | 4.067720E+2 | 32.62160788, -1.47156E+2, -0.70775567, 1.32217E+15, -1.98453E+1 | 4.0 |
| Antigone | 33.82722561 | 6.84665005, 9.899551E+5, 2.71799837, -9.89955E+5, 2.71799837 | 64.0 |
| Couenne | 4.508885E+2 | 56.57155475, 2.71247E-10, 0.11751213, -8.03590E+1, -0.22170130 | 8.0 |
| Lingo | 25.25308277 | 5.17113395, 2.36035084, 2.40494700, -2.34875359, 2.40601371 | 3.0 |
| Scip | 25.25308388 | 5.17112097, -3.25522448, 2.40586326, 3.26682183, 2.40509306 | 1.0 |
| Matlab-GA | 411.5534198 | 38.60769300, -118.39477447, -0.51680714, 540.13847085, -.18503202 | 7.5 |
| 1stOpt | 0.909227034 | -1.49546149, 1.35117337, 0.76999609, -1.28596489E-13, 7.27408071 | 7.0 |

Appx. Table 15: Detailed Results of Problem 15

| Solver | Obj. Values | 6 Parameters | Time (s) |
| --- | --- | --- | --- |
| Baron | 0.29936497 | 89.64131684, 0.00034900, 0.08923089, 1.152045E+2, 45.90368506, -0.03534606 | 1.0 |
| Antigone | 1.75498593 | 86.47537694, -0.06079708, 0.02388325, 8.615526E+7, -2.29546E+8, -0.06184580 | 5.0 |
| Couenne | 9.607825E+4 | 0.00000000, 0.00000000, 0.00000000, 0.00000000, 0.00000000, 0.00000000 | |

| Lingo | 0.98870723 | 19.14763322, -0.50000142, 0.00334531, 53.58551771, 0.00004210, 0.11457453 | 11.0 |
| --- | --- | --- | --- |
| Scip | 6.586250E+2 | 1.811573E+4, 1.397261E+4, -8.30737E+2, -1.80225E+4, -1.99696E+2, -1.94744E+3 | 1.0 |
| Matlab-GA | 68.55265925 | 87.95324350, 0.05462706, 0.03971266, 59.05571081, 6.57603868, -0.04870532 | 5.5 |
| 1stOpt | 0.299364969 | 89.64131554, 0.00034900, 0.08923089, 115.20457437, 45.90369795, -0.03534605 | 3.0 |

Appx. Table 16: Detailed Results of Problem 16

| Solver | Obj. Values | 6 Parameters | Time (s) |
| --- | --- | --- | --- |
| Baron | 0.03495035 | 3.66108963, 0.12571991, -0.00010718, -0.00000023, -7.60002442, 0.61030587 | 2.5 |
| Antigone | 0.06097539 | 3.37886196, -7429.15, -0.0026876, 7429.357, -0.99999276, 0.63714881 | 3.0 |
| Couenne | 0.09270095 | 4.18124473, 0.0750258, -0.369909, -0.03390706, 1.83138218, -1.47349466 | 1.5 |
| Lingo | 0.06097754 | 3.37877681, -379.688, -0.00267922, 379.8948, -0.99985833, 0.63719085 | 1.0 |
| Scip | 7.07618234 | 348.1798, -0.01766722, 0.10522572, -4.22029998, 0.08294398, 0.85904795 | 10.0 |
| Matlab-GA | 0.28979101 | 6.89667282, 0.08000000, 0.34807958, -0.37989770, -0.10160383, 1.22222846 | 3.6 |
| 1stOpt | 0.00712837 | -1.10289481E-7, -0.06291472, 0.01039641, 21.45616947, 0.04512023, 217.96391230 | 8.0 |

Appx. Table 17: Detailed Results of Problem 17

| Solver | Obj. Values | 5 Parameters | Time (s) |
| --- | --- | --- | --- |
| Baron | 4.301794E+4 | 1.00000E+10, 0.85333953, -0.14580455, 1.662712E+2, -3.34923342 | 23.0 |
| Antigone | 7.555325E+7 | 0.00000000, 0.00000000, 0.00000000, 0.00000000, 0.00000000 | 1 |
| Couenne | 1.694418E+7 | 3.194804E+3, 17.10800000, 1.940918E+8, -2.10000000, 1.00000000 | 1.0 |
| Lingo | 7.974938E+4 | 7.780503E+7, -2.50000E+9, -2.44338E+9, 6.925638E+2, -1.84715161 | 15.0 |
| Scip | 7.469755E+7 | 0.00019489, 1.00316225, 0.00317225, -1.03510E+1, 1.00000175 | 1.0 |
| Matlab-GA | 158654.0395 | 382.23268780, 663.73749754, 15.24850232, -3.13400036, 0.87681539 | 7.1 |
| 1stOpt | 34105.93604 | 978429282056.723, 0.85603224, -0.148796, 101.1984159, -4.71601407 | 3.5 |

Appx. Table 18: Detailed Results of Problem 18

| Solver | Obj. Values | 6 Parameters | Time (s) |
| --- | --- | --- | --- |
| Baron | 0.00551218 | -1.14314E+2, -3.80969028, 65.50591494, -0.00138626, 0.13612078, -6.53153E+1 | 1.5 |
| Antigone | 0.00693548 | -9.65760E+1, -3.65786062, -4.28675E+5, 0.00000002, -0.00002268, 4.286747E+5 | 5.0 |
| Couenne | 0.00885065 | -1.71647E+2, -4.46163263, 14.48185594, -1.04407097, 0.55786212, 0.22946538 | 0.5 |
| Lingo | 0.01876434 | 29.66746159, -0.81018709, 61.21006510, -0.03166914, -0.23077972, -5.32738E+1 | 0.3 |
| Scip | 0.13393480 | -1.84009E+2, -0.05786949, 5.385833E+5, -0.00007248, -0.00001633, -5.38355E+5 | 2.0 |
| Matlab-GA | 0.11187831 | -39.52428458, -7.27143440, 137.01074921, -2.13101494, -0.94217367, 0.65051835 | 1.9 |

| Solver | | | |
|---|---|---|---|
| 1stOpt | 0.00551218 | -114.31450705, -3.80969036, 65.50588759, -0.00138626, 0.13612083, -65.31525149 | 4.0 |

Appx. Table 19: Detailed Results of Problem 19

| Solver | Obj. Values | 5 Parameters | Time (s) |
|---|---|---|---|
| Baron | 0.0012525071 | 2.800484, -0.252724, 0.01714066, -1.1677863470, 0.5195859673 | 1.5 |
| Antigone | 0.01661039 | 3.889934E+4, 1.472099E+5, -2.39497E+4, 1.78739419, -0.11036775 | 1.0 |
| Couenne | 0.26418161 | -0.00128909, -0.40676589, 0.04084441, 2.80579502, -0.26006925 | 1.0 |
| Lingo | 0.012957474 | 0.18830221, 0.36225587, -0.8601513E-01, 1.945311507, -.129558395 | 1.0 |
| Scip | 0.26418161 | -0.00128909, -0.40676589, 0.04084441, 2.80579502, -0.26006925 | 1.0 |
| Matlab-GA | 0.02702564 | 17.29731664, 72.62639522, -11.68396586, 1.63329362, -0.07439744 | 5.0 |
| 1stOpt | 0.0012525071 (80% Pr.) | 2.800484, -0.252724, 0.01714066, -1.1677863470, 0.5195859673 | 2.5 |

Appx. Table 20: Detailed Results of Problem 20

| Solver | Obj. Values | 4 Parameters | Time (s) |
|---|---|---|---|
| Baron | 0.00086491 | 2.45379460, 0.00000000, 1.88303559, 5.26129573 | 1.0 |
| Antigone | 0.56099429 | 3.121775E+5, -1.14161888, -0.00000076, -9.77670E+5 | 1.0 |
| Couenne | 0.35950147 | -1.16114E+1, 1.057125E+2, -0.28919925, 8.21986413 | 1.0 |
| Lingo | 0.05563105 | -2.98702836, 39062501.2297, -1.2622715, 7.66175799 | 70.0 |
| Scip | 9.56709091 | -3.12494E+5, -0.00000009, -4.95150278, 3.125008E+5 | 8.0 |
| Matlab-GA | 0.41550289 | -23.30781084, 51.17599859, -0.21298225, 8.15429463 | 3.0 |
| 1stOpt | 0.00086490 (80% Pr.) | 2.45379459, 2.61800799E-12, 1.88303559, 5.26129573 | 2.0 |

Appx. Table 21: Detailed Results of Problem 21

| Solver | Obj. Values | 6 Parameters | Time (s) |
|---|---|---|---|
| Baron | 0.18149434 | 1.71260045, 2.727011E+5, -0.56767429, -6.93866083, 15.00250094, -3.26254261 | 78.8 |
| Antigone | 4.16891524 | -2.09698E+5, 2.097201E+5, -0.00164362, 1.01538792, 0.02767342, 1.189979E+4 | 1.0 |
| Couenne | 33.62486495 | 6.90367883, 5.00679E-11, -0.32450746, 0.999905E+5-2.93270E+1, -1.59994E+1 | 1.0 |
| Lingo | 0.21067279 | 2.17931926, 9.190488E+3, -0.27546037, -3.11628838, 17.13682655, -4.21057831 | 20.0 |
| Scip | 26.91510583 | 7.34510077, 0.24721349, -0.82364514, -2.08982E+1, -3.50248732, -6.35319E+1 | 3.0 |
| Matlab-GA | 7.73687752 | 3.00939311, 168.72279853, -0.92629491, -55.02290357, 3.40948175, -10.92574611 | 10.0 |
| 1stOpt | 0.181494338 (80% Pr.) | 1.71260230, 272701.205533, -0.56767387, -6.93865237, 15.00250991, -3.26254613 | 5.0 |

Appx. Table 22: Detailed Results of Problem 22

| Solver | Obj. Values | 6 Parameters | Time (s) |
|---|---|---|---|
| Baron | 3.399997E+3 | -2.08932E+2, 1.95000105, -2.94044740, 0.01315072, -2.16340504, -0.00667851 | 500.0 |
| Antigone | 5.885228E+3 | -1.04390E+2, 0.91123275, 6.875000E+5, -1.00000E+6, -1.50488593, -0.00350768 | 25.0 |

| Couenne | 5.567073E+3 | -0.00000144, -1.12122E+1, 4.656477E+2, 0.01514596, 0.52280678, 0.01838243 | 2.5 |
|---|---|---|---|
| Lingo | 483.32350 | 336.50619, -0.372915E-01, 0.10000E+11, -0.15517375, 41054012.932, -0.15163727 | 1.0 |
| Scip | 8.108606E+5 | -4.52314E+2, 6.31871343, -0.68985984, 0.04403902, -0.00264650, 0.04438553 | 2.0 |
| Matlab-GA | 8378.444731 | 149.44676959, 1.25128791, -10.15528486, 0.010, -742.96027136, -0.36251651 | 3.0 |
| 1stOpt | 179.3335177 (50% Pr.) | 235.74202465, -0.58139127, 6.956965E-10, 0.23269268, 2.14678E-12, 0.23265419 | 4.0 |

Appx. Table 23: Detailed Results of Problem 23

| Solver | Obj. Values | 6 Parameters | Time (s) |
|---|---|---|---|
| Baron | 2.12289246 | 7.662886E+7, 9.919645E+9, 1.85419863, 9.440589E+5, 0.00000093, 6.00359081 | 1.0 |
| Antigone | 5.68622712 | -2.02833E+3, -0.59137709, 0.29308656, 0.04539261, 3.476075E+3, -0.38430063 | 0.5 |
| Couenne | 8.14877184 | 34.16682626, 5.113982E+4, 2.51760775, 0.31899285, 0.03985195, 1.65561449 | 3.0 |
| Lingo | 2.12333274 | 6.79914E+11, 9.08728E+13, 1.86198513, 8.378249E+9, 0.00000094, 5.99853452 | 1.0 |
| Scip | 1.132110E+3 | 0.00000000, 0.00000000, 0.50000000, 0.00000000, 0.00000000, -1.44269504 | 1.0 |
| Matlab-GA | 17.86916266 | 32.22348125, 240.48294814, 1.24511497, 0.20509204, -2.94714343, -0.72351801 | 17.8 |
| 1stOpt | 1.978953145 (90% Pr.) | -50.76268102, -3665.045429, 1.67897488, -0.78081813, 4.6179945E-5, 4.53154207 | 4.0 |

Appx. Table 24: Detailed Results of Problem 24

| Solver | Obj. Values | 5 Parameters | Time (s) |
|---|---|---|---|
| Baron | 0.99685339 | 1.828909E+3, 4.694089E+2, -4.75041E+2, 0.99857040, -8.46839E+5 | 1.0 |
| Antigone | Failed | | |
| Couenne | 0.99685339 | 1.828911E+3, 4.694090E+2, -4.75041E+2, 0.99857040, -8.46840E+5 | 1.0 |
| Lingo | 3.6361175 | 283.07944, 0.3768249E-01, 6574443.0105, -2.330047145, -32909.6022 | 1.0 |
| Scip | 3.936153E+5 | 13.14988661, -0.01841645, -2.20260279, -0.21227055, 0.01309168 | 2.0 |
| Matlab-GA | 40.36491692 | -71.13251807, -0.02761794, 80.2134613, 0.23131496, -7863.4636352 | 31.0 |
| 1stOpt | 0.52007113 (20% Pr.) | -372.438433, -1989.2957928, 1992.77996, 0.99974684, -651705.52987 | 11.0 |

Appx. Table 25: Detailed Results of Problem 25

| Solver | Obj. Values | 6 Parameters | Time (s) |
|---|---|---|---|
| Baron | 0.44499232 | 0.86085780, 26.39091755, -4.20854667, 0.00000000, 0.00031008, -0.37943320 | 1.0 |
| Antigone | 2.79312763 | -4.66112129, -3.64374E+2, -1.82325586, -0.00016337, -1.24702E+1, 1.938834E+2 | 0.5 |
| Couenne | 2.79312763 | -4.66112129, -2.44732E+1, -1.82325586, -0.00001097, -1.77039948, 9.94725960 | 1.5 |
| Lingo | 1.43570424 | 28.13228333, 1.608090E+8, 1.01507793, -4.99851E+9, -0.34938283, -0.39126432 | 118.0 |
| Scip | 7.891489E+2 | 18.68888889, -1.12410E+2, -7.92507378, 2.76111045, 55.59497440, 5.20803340 | 195.0 |

| Solver | Obj. Values | 6 Parameters | Time (s) |
|---|---|---|---|
| Matlab-GA | 6.80895741 | -20.67085107, 881.07869097, -0.82335712, 0.08683143, 11.18127667, 0.30365670 | 14.0 |
| 1stOpt | 0.444992322 (50% Pr.) | 0.86085782, 26.39091755, -4.2085467, -4.6205066E-11, 0.0003100812, -0.3794332 | 8.0 |

Appx. Table 26: Detailed Results of Problem 26

| Solver | Obj. Values | 6 Parameters | Time (s) |
|---|---|---|---|
| Baron | 0.08897309 | -3.67039077, -2.67616E+2, 40.14381031, 0.00000068, 28.59370576, 0.00298226 | 1.0 |
| Antigone | 0.62251302 | 0.07478335, 9.804494E+9, 61.91498574, 0.00331583, 14.99312007, 0.00014198 | 475.0 |
| Couenne | 0.08897309 | -3.67039077, -2.67616E+2, 40.14381031, 0.00000068, 28.59370577, 0.00298226 | 1.0 |
| Lingo | 0.62251342 | 0.09062279, 1.758288E+8, 61.91519251, 0.00331584, 14.99313690, 0.00014198 | 1.0 |
| Scip | 1.35478951 | -7.50582351, -2.62156E+2, -1.07589E+5, 0.00002368, 1.076548E+5, 0.00002365 | 1.0 |
| Matlab-GA | 78.40704417 | 0.78831320, 10086.4714999, 58.6185013, 0.00191607, 237.01935085, 189.23241894 | 25.0 |
| 1stOpt | 0.005260412 (90% Pr.) | -13.8111306, -191.2019687, 127.0998036, -4.322261E-6, -26.3906786, 0.00039558 | 3.0 |

Appx. Table 27: Detailed Results of Problem 27

| Solver | Obj. Values | 5 Parameters | Time (s) |
|---|---|---|---|
| Baron | 3.289399E+3 | 1.825534E+2, -2.07516E+6, -0.00102456, -3.95933E+3, -2.00094E+3 | 3.0 |
| Antigone | 3.291632E+3 | 54.15333332, 48.40927591, -0.03961487, 3.633255E+2, 2.760160E+2 | 65.0 |
| Couenne | 4.317044E+3 | 58.53558240, -4.09404E+3, -0.01695870, -2.71499E+1, 12.87420466 | 1.5 |
| Lingo | 4.375858E+3 | 1.136729E+2, 86.84000000, -0.00367159, -1.0000E+10, 1.688385E+2 | 793.0 |
| Scip | 7.646795E+3 | 57.40940455, 0.01673885, 0.00000000, 0.00000000, 0.00000000 | 538.0 |
| Matlab-GA | 4433.438896 | 56.09856009, 41.78568256, 0.04407766, 0.25901240, 51.37769758 | 2.0 |
| 1stOpt | 6.681978318 (40% Pr.) | 127.98329215, -394.85726762, -0.01292565, 177.29323334, 200.39679 | 3.5 |

Appx. Table 28: Detailed Results of Problem 28

| Solver | Obj. Values | 6 Parameters | Time (s) |
|---|---|---|---|
| Baron | 0.65965216 | -7.54073E+3, -0.00250305, 2.030577E+3, 7.564940E+3, -0.00249222, 2.030578E+3 | 196 |
| Antigone | No solution returned | | |
| Couenne | 0.02493958 | -0.00480263, 0.01755207, 2.034109E+3, 24.30667862, 0.00039440, 2.028803E+3 | 8.0 |
| Lingo | 2.84091799 | 2.11251E-14, 0.00001241, 4.034780E+2, 1.979110E+4, -0.00000107, -5.04235E+2 | 1.0 |
| Scip | 1.090121E+4 | 0.00000000, 0.00000000, 1.275005E+3, 0.00000000, 0.00000000, 1.275005E+3 | 1.0 |
| Matlab-GA | 3.65573320 | 24.88769031, -0.00003205, 2044.0331364, 2995.169635, -7.22825095, 2011.06010296 | 8.0 |
| 1stOpt | 0.02493958 (40% Pr.) | 24.30668467, 0.00039442, 2028.80292028, -0.00480442, 0.01754741, 2034.11118702 | 5.0 |

Appx. Table 29: Detailed Results of Problem 29

| Solver | Obj. Values | 5 Parameters | Time (s) |
|---|---|---|---|
| Baron | 0.00073779 | -0.56850923, -1.83104E+1, 15.23082574, 3.600481E+6, 0.00756081 | 6.0 |
| Antigone | 0.70660012 | -1.45986E+1, 0.63935818, -1.61981E+1, 2.94864376, 0.00000001 | 1.0 |
| Couenne | 0.72451259 | -2.55275140, 1.000020E+5, -5.49295501, 0.12092048, 0.00814697 | 6.5 |
| Lingo | 0.49538789 | -6.07359041, -0.87642035E-01, -4.49800535, 68.28569489, -0.42198461 | 1.0 |
| Scip | 0.72451259 | -2.55275161, 6.013853E+6, 0.76182399, 62.93844641, 0.00814696] | 4.0 |
| Matlab-GA | 0.38729520 | 1.00265269, 62.52923231, 15.18933746, 505.57673949, 0.49233857 | 4.5 |
| 1stOpt | 0.000737787 (70% Pr.) | -0.56850925, -18.310435, 15.2308257, 3600482.61782988, 0.0075608 | 2.0 |

Appx. Table 30: Detailed Results of Problem 30

| Solver | Obj. Values | 6 Parameters | Time (s) |
|---|---|---|---|
| Baron | 0.66930130 | 1.305701E+3, -2.17482E+3, -0.00926567, 0.02059264, -2.83409E+1, -0.04629415 | 96.0 |
| Antigone | 0.46183999 | -9.93543E+1, 98.51437190, -0.00105731, -0.00192505, -1.06147898, -0.00194778 | 1.0 |
| Couenne | 0.46183999 | -9.93543E+1, 98.51437190, -0.00105731, -0.00192505, -1.06147898, -0.00194778 | 3.0 |
| Lingo | 0.651761430 | -1697.47777663, 1697.47454408, -0.25462635E-06, -0.1232717E-03, -1.00024846218, -0.12362196E-03 | 1.0 |
| Scip | 2.476439E+4 | 1.206837E+2, -1.73802E+1, 0.00977466, 0.00125971, -0.18990088, 0.00978600 | 1.0 |
| Matlab-GA | 999.2931747 | 355.7818135, 7165.54968299, -0.05690627, 0.00472166, 68.59190142, -0.06165643 | 20.0 |
| 1stOpt | 0.461594701 (50% Pr.) | 75.0385145, -75.9037924, 0.00144625, 0.00209540, -1.06385958, 0.00190959 | 3.0 |

Appx. Table 31: Detailed Results of Problem 31

| Solver | Obj. Values | 7 Parameters | Time (s) |
|---|---|---|---|
| Baron | 2.24560393 | 1.538366E+5, 8.41349914, 0.01135526, 10.53351816, 6.18101967, -9.29640876, -0.42010569 | 2.0 |
| Antigone | 5.38894747 | 0.00015193, -5.54886130, 0.00624608, 63.79909007, 15.95843724, -8.48406E+1, 0.10987780 | 1.0 |
| Couenne | 8.214714E+4 | -6.1586E+13, 5.62143069, -1.0979E+11, 0.79909098, 0.00012542, 1.000000E+5, 0.999546E+5 | 39.0 |
| Lingo | 13.93232765 | 0.00123035, -2.58904089, -0.64781715, 8.16246051, 0.03152509, -8.61293229, 0.00196161 | 1.0 |
| Scip | 1.000686E+5 | -1.37918E+2, -4.22329318, -4.42260E+4, 1.688906E+6, 2.60506147, -1.68890E+6, 0.00000055 | 30.0 |
| Matlab-GA | 17.07814305 | 49.96974468, -0.81912836, 0.46127996, 2.00056237, -8.08924422, 2.63588091, 3.07030468 | 6.5 |
| 1stOpt | 2.24560393 (50% Pr.) | 153845.38319049, 8.41353246, 0.01135526, 10.5336822, 6.18102749, -9.29652312, -0.42010346 | 15.0 |

Appx. Table 32: Detailed Results of Problem 32

| Solver | Obj. Values | 5 Parameters | Time (s) |
|---|---|---|---|
| Baron | 1.742015E+2 | -2.34215503, 0.00080581, 29.78653339, 0.51331490, -0.19256104 | 1.0 |
| Antigone | 1.742015E+2 | -2.34215503, 0.00080581, 29.78653339, 0.51331490, -0.19256104 | 1.0 |
| Couenne | 1.742015E+2 | -2.34215503, 0.00080581, 29.78653339, 0.51331490, -0.19256104 | 1.0 |
| Lingo | 1.742015E+2 | -2.34215503, 0.00080581, 29.78653339, 0.51331490, -0.19256104 | 1.0 |

| Solver | | | |
|---|---|---|---|
| Scip | 7.084184E+3 | -2.73251902, 5.23685237, -0.05606639, 3.52535022, 0.01195560] | 600.0 |
| Matlab-GA | 183.6248605 | 6.07252678, 37.56733705, 32.08049779, 0.53528178, -0.20725961 | 4.0 |
| 1stOpt | 8.31641750 (40% Pr.) | -5.60698775, -0.00106695, -178.65762017, -0.19318011, -0.06534412 | 3.0 |

Appx. Table 33: Detailed Results of Problem 33

| Solver | Obj. Values | 5 Parameters | Time (s) |
|---|---|---|---|
| Baron | 2.46322E-20 | 17.75258710, 2.753477E+7, -3.45700E+1, 1.25216863, 0.26178963 | 4.0 |
| Antigone | 6.89358250 | -0.75763265, -0.00000743, -2.13292732, -0.00000001, -3.98339E+4 | 1.0 |
| Couenne | 39.02123764 | 7.62499962, 1.99999987, 0.00394471, 0.76544553, 1.88998692 | 1.0 |
| Lingo | 0.84093052 | 2.18344917, 4.86551975, -1.61540E+2, -0.01250211, 0.70663750 | 1.0 |
| Scip | 1.090936E+5 | 2.24607095, 4.43028516, -0.73607442, 0.01190770, 0.49106554 | 1.0 |
| Matlab-GA | 8.64425832 | 30.09874849, 46.63058508, -19.53323872, 28.40063804, 0.72000000] | 10.0 |
| 1stOpt | 5.16987E-26 (90% Pr.) | 17.752587, 27534771.2015618, -34.57000754, 1.25216863, 0.26178963 | 10.0 |

Appx. Table 34: Detailed Results of Problem 34

| Solver | Obj. Values | 5 Parameters | Time (s) |
|---|---|---|---|
| Baron | 0.01228654 | 22.75972413, 14.51011916, 10.62880288, 0.06203956, -4.37482808 | 1.0 |
| Antigone | 0.14654921 | 7.73933133, -0.01674006, 23.24398802, -0.13755821, 0.50872649 | 1.0 |
| Couenne | 0.01237477 | 10.39037779, 0.07561919, 22.75812030, -0.06168785, -0.18340909 | 1.0 |
| Lingo | 0.01237477 | 10.39037779, 0.07561919, 22.75812030, -0.06168785, -0.18340909 | 1.0 |
| Scip | 0.01228654 | 22.75972406, 14.51011952, 10.62880289, 0.06203956, -4.37482889 | 1.0 |
| Matlab-GA | 0.35758658 | 32.62872681, 3.830000, 10.47075253, 0.050000,43.64959327 | 1.0 |
| 1stOpt | 0.010808142 (50% Pr.) | 514.22809651, 0.00131837, 23.26974837, 0.05947556, 38.57101985 | 3.0 |

Appx. Table 35: Detailed Results of Problem 35

| Solver | Obj. Values | 5 Parameters | Time (s) |
|---|---|---|---|
| Baron | 75.57370284 | 95.95391750, -2.67239E+8, 1.688022E+2, -1.69188E+2, 0.99965710 | 53.0 |
| Antigone | 3.296046E+2 | -1.49896E+6, 1.498393E+6, -0.06329470, 0.06330491, 0.99997574 | 1.0 |
| Couenne | 4.625350E+2 | -1.95733E+1, 5.95637121, 0.01487598, -0.00000006, 2.99436779 | 20.0 |
| Lingo | 329.78055817 | -119982.3155, 119416.3498, -17.22854669, 17.2286742, 0.99999888 | 1.0 |
| Scip | 6.570858E+3 | 62.19454545, 0.18917211, -1.85395263, -0.98967486, -2.01314235 | 8.0 |
| Matlab-GA | 814.3733233 | -152.296004, 5849.316205, -0.00603119, -671.0336105, -1.08343279 | 35.0 |
| 1stOpt | 30.02585461 (60% Pr.) | 406.470778, -406.674725, -0.00085920, 7.5085119E-18, 6.33124939 | 6.0 |

Appx. Table 36: Detailed Results of Problem 36

| Solver | Obj. Values | 5 Parameters | Time (s) |
|---|---|---|---|
| Baron | 0.81336718 | 7.142952E+8, -0.49631409, -3.15526E+2, 0.05258203, 3.98232762 | 1.0 |
| Antigone | 1.461242E+3 | 4.969602E+5, 3.18727781, -2.31686E+2, -0.00000619, -4.98119E+5 | 2.0 |
| Couenne | 1.178249E+3 | 1.43559987, -7.00285E+2, -3.28599E+2, -0.01632977, 25.12638863 | 1.0 |
| Lingo | 1330.351267 | -469.600776, -184.2795488, -292.726348, 0.4533577E-02, 114.313141 | 1.0 |
| Scip | 1.085889E+4 | -0.24921566, 0.00008423, 0.00000000, 7.22546673, 46.65100000 | 1.0 |

| Matlab-GA | 1606.349399 | 1059.0105589, 1607.7431426, -261.1912764, 0.01072965, -21.455213 | 23.5 |
| 1stOpt | 0.81336718 (80% Pr.) | 714295206.87674, -0.49631409, -315.526295, 0.05258203, 3.9823278 | 2.0 |

Appx. Table 37: Detailed Results of Problem 37

| Solver | Obj. Values | 5 Parameters | Time (s) |
| --- | --- | --- | --- |
| Baron | 0.00085308 | -1.06145E+2, 1.387227E+2, -0.00822761, 1.408015E+2, -1.12121E+2 | 2.0 |
| Antigone | 0.00069673 | -2.35387250, 1.43051342, 0.00841358, -6.31155E+1, -5.35261E+2 | 150.0 |
| Couenne | 0.00069673 | -2.35387251, 1.43051342, 0.00841358, -6.31155E+1, -5.35261E+2 | 12.0 |
| Lingo | 0.00388164 | -2.73544373, 8.58092576, -0.00037923, -3.72008227, 26.7422079 | 0.5 |
| Scip | 0.00085308 | -1.06145E+2, 1.387227E+2, -0.00822761, 1.408015E+2, -1.12121E+2 | 1.0 |
| Matlab-GA | 0.00868975 | 1.64263749, -96.67572976, -0.42420802, -0.42420802, 85.39951546 | 1.5 |
| 1stOpt | 0.00069673 (70% Pr.) | -2.35387248, 1.43051340, 0.00841358, -63.11547564, -535.26137948 | 2.5 |

Appx. Table 38: Detailed Results of Problem 38

| Solver | Obj. Values | 6 Parameters | Time (s) |
| --- | --- | --- | --- |
| Baron | 7.176911E+2 | -1.41886E+2, -2.66199E+3, -0.02095992, -2.50789E+2, -0.00003080 | 1.0 |
| Antigone | 7.176883E+2 | -1.41896E+2, -4.07721E+6, -3.20581E+1, -3.82455E+5, -0.00003084 | 1.0 |
| Couenne | 1.475831E+6 | 5.617500E+2, 1.29293241, 1.569623E+3, 41.53826697, 0.00000000 | 1.0 |
| Lingo | 7.165598E+2 | -1.44770E+2, 16.21236175, 0.00008265, 0.33974452, -0.00004931 | 1.0 |
| Scip | 1.584736E+6 | 4.450752E+2, 0.00000000, 0.00000000, 0.00000000, 0.00004076 | 1.0 |
| Matlab-GA | 1199841.698 | 1239.62628, 11371.924, 6.10731959e-04, 107.16827289, 105.21283 | 23.0 |
| 1stOpt | 676.71429733 (60% Pr.) | -239.493475, 48.425682, 3.2478409E-5, 5.4879009E-19, -0.00123899 | 2.3 |